\documentclass[11pt]{article}

\usepackage{makeidx}
\usepackage{amssymb}
\usepackage{amsfonts}
\usepackage{amsmath}
\usepackage{euscript}
\usepackage{theorem}
\usepackage{enumerate}
\usepackage{graphics}

\theorembodyfont{\rm}      

\theoremstyle{change}      

\begingroup

\newtheorem{thm}{Theorem\hskip 5mm}[section]

\newtheorem{prop}[thm]{Proposition\hskip 5mm}

\newtheorem{cor}[thm]{Corollary\hskip 5mm}

\newtheorem{lem}[thm]{Lemma\hskip 5mm}

\newtheorem{exa}[thm]{Example\hskip 5mm}

\newtheorem{note}[thm]{Note\hskip 5mm}

\endgroup

\begingroup

\endgroup

\begingroup

\endgroup

\def\Sp{{\rm Sp}}
\def\GL{{\rm GL}}

\def\l{{\lambda}}

\def\p{{\phi}}
\def\e{{\ell}}

\def\P{{\cal P}}
\def\HP{{\check{{\cal P}}}}

\def\t{{\vartheta}}
\def\u{{\varepsilon}}

\begin{document}

\begin{center}
{\bf Steinberg lattice of the general
linear group and its modular reduction}
\end{center}

\begin{center}
{\sc Fernando Szechtman}
\end{center}

\section{Introduction}

Let $G=\GL_n(q)$ be the general linear group of degree $n\geq 2$
defined over a finite field $F_q$ of characteristic $p$. We fix a
prime $\ell\neq p$ and let $R$ denote a local principal ideal
domain having characteristic 0, maximal ideal $\ell R$, and
containing a primitive $p$-th root of unity. Then the residue
field $K=R/\ell R$ has characteristic $\ell$ and a primitive
$p$-th root of unity.

By a Steinberg lattice of $G$ over $R$ we understand a left
$RG$-module, say $M$, which is free of rank $q^{n(n-1)/2}$ as an
$R$-module and affords the Steinberg character. The reduction of
$M$ modulo $\ell$ is the $KG$-module $M/\ell M$.

In this paper {\em the} Steinberg lattice is the left ideal $I=RG\cdot
\u$ of the group algebra $RG$,
\begin{equation}
\label{defechico} \u=\underset{\sigma \in
S_n}\sum\mathrm{sg}(\sigma)\sigma\widehat{B},
\end{equation}
where the symmetric group $S_n$ is viewed as a subgroup of $G$, we let $B$
denote the upper triangular group, and
$\widehat{S}=\underset{s\in S}\sum s$ for any subset $S$ of $G$.

Our main object of study is the
the $\ell$-modular reduction of $I$, namely the $KG$-module $L=I/\ell I$.
In particular, we wish to find a composition series of $L$,
the socle and radical series of $L$, the length, say $c(L)$, of $L$, and any additional
structural information about  $I$ and $L$ that might be of use in
achieving these goals, or interesting in its own right.

Many other Steinberg lattices and their corresponding reductions modulo $\ell$ appear in a natural manner,
and will be compared to $I$ and $L$.

The first results are due to Steinberg [5]. Let $U$ be the upper unitriangular group, i.e. the Sylow
$p$-subgroup of $B$. Then $I$ is a free $R$-module with basis $\{u\u\,|\, u\in U\}$ and $U$ acts on $I$ via the
regular representation. Naturally $L$ has $K$-basis $\{u\cdot (\u+\ell I)\,|\, u\in
U\}$ and affords the regular representation of $U$. Moreover, $L$ is irreducible if and only if $\ell\nmid [G:B]$.
Steinberg did not state it explicitly, but it is obvious from [5] that the socle of $L$, say $soc(L)$, is irreducible.

There is a canonical symmetric bilinear form $RG\times RG\to R$
given by $(g,h)\mapsto \delta_{g,h}$. Restriction to $I$ followed
by scaling by $1/|B|$ yields the $G$-invariant symmetric bilinear
form $f:I\times I\to R$ with zero radical studied by Gow in [4]. He uses
$f$ to produce the $RG$-submodules $I(c)$ of $I$ given by
$$
I(c)=\{x\in I\,|\, f(x,I)\subseteq \ell^c R\},\qquad c\geq 0.
$$
This yields the following filtration of $RG$-modules, where all inclusions are strict:
\begin{equation}
\label{pato1}
I=I(0)\supset I(1)\supset I(2)\supset\cdots.
\end{equation}
He next considers the $KG$-submodules $L(c)$ of $L$ defined by
$$
L(c)=(I(c)+\ell I)/\ell I,\quad c\geq 0,
$$
which give rise to a filtration of $KG$-modules
\begin{equation}
\label{pato2}
L=L(0)\supseteq L(1)\supseteq L(2)\supseteq\cdots.
\end{equation}
Each factor of (\ref{pato2}) is a $KG$-module and will be denoted by
$$M(c)=L(c)/L(c+1),\quad c\geq 0.
$$

As $L$ is finite dimensional, the series (\ref{pato2}) must
eventually stabilize and there may be prior repetitions. The
question as to when exactly this happens was settled by Gow.
Write $\P$ for the lattice of standard
parabolic subgroups of $G$, i.e. those containing~$B$. A non-negative integer
$c$ is said to be a $\P$-value if $\ell^c\mid [G:P]$ but $\ell^{c+1}\nmid [G:P]$, i.e.
$\nu_\ell([G:P])=c$, for some $P\in\P$. Let $V$ stand for the total number of $\P$-values.
Gow proves that the factor $M(c)$
is non-zero if and only if $c$ is a $\P$-value. Furthermore, if
$b=\nu_\ell([G:P])$, the largest $\P$-value, then $L(b)=soc(L)$, which by above is irreducible.
Since $L(b)/L(b+1)\neq 0$, it follows that $0=L(b+1)=L(b+2)=\cdots$.

Clearly Gow's work implies $c(L)\geq V$, with equality if and only
if $M(c)$ irreducible for every $\P$-value  $c$, that is, if and
only if (\ref{pato2}) is a composition series of $L$. All repeated
terms in (\ref{pato2}) must be deleted when interpreting this
statement.  Looking at the last line of the decomposition matrices
for unipotent representations of $\GL(n,q)$, $n\leq 10$, as given
by James in [3], Gow believed that $c(L)=V$ and conjectured this
would hold for any $n$.

Let
$$
e=\mathit{min}\{\,i\geq 2\,|\, \ell\text{ divides }\frac{q^i-1}{q-1}\},
$$
and note that if
$\ell\nmid q-1$ then $e$ divides $\ell-1$ and is the order of
$q$ modulo $\ell$,  while if $\ell\mid q-1$ then $e=\ell$.

Gow's observation is based on the matrices explicitly displayed in [3], which equal the
decomposition matrices as long as $\lfloor n/e\rfloor<\ell$. To obtain the
latter when $\lfloor n/e\rfloor\geq \ell$ requires adjustment matrices,
as indicated by James. We will come back to this point shortly.

As long as $\lfloor n/e\rfloor<\ell$, Gow's conjecture does hold for any $n$, as shown by Ackermann (see Section 4.6 of [1]), who proved, among
many other things, that $L$ is uniserial of length $c(L)=V=\lfloor n/e\rfloor+1$, provided $\lfloor n/e\rfloor<\ell$.
Furthermore, Theorem \ref{A1} verifies Gow's conjecture in many other cases, including the case
$\lfloor n/e\rfloor\leq \ell$, and Theorem \ref{A2} proves that if $\lfloor n/e\rfloor>\ell$
then the first $\ell+1$ non-zero factors of (\ref{pato2}), starting at the bottom, are indeed irreducible.
In addition, Theorem \ref{A3} proves that the top factor $M(0)=L(0)/L(1)$ is irreducible under no restrictions at all.
Moreover, we know from [6] that $M(c)$ is a completely reducible $KG$-module. Furthermore,
Sections \ref{pre} and \ref{nulam} associate a non-zero cyclic submodule $N(P)$ of $M(c)$ to any
$P\in\P(c)$ and prove it to be irreducible.

In spite of all this evidence Gow's conjecture is actually false.
Indeed, for $c\geq 0$ let $\P(c)$ consist of all $P\in\P$ such
that $\nu_\ell([G:P])=c$. Let $\P^*$ be the set of all standard
parabolic subgroups that correspond to partitions of $n$ where
every part is either 1 or of the form $e\ell^i$ for some $i\geq
0$. Define $\P^*(c)=\P(c)\cap\P^*$. Section \ref{cmp} shows that
$M(c)$ equals the direct sum of all distinct $N(P)$ as $P$ runs
through $\P^*(c)$. Thus, Gow's conjecture transates as follows:
the $N(P)$ are all {\em equal}, $P\in\P^*(c)$, whenever $c$ is a
$\P$-value. We tried very hard to prove this, without success.
Then, together with D. Djokovic, we examined James' tables and
noticed that the correct decomposition matrices for $n\leq 10$
give $c(L)=\vert\P^*\vert$, which is equivalent to the $N(P)$
being {\em distinct} for all $P\in\P(c)$ and all $\P$-values $c$,
i.e. a composition series of $L$ is obtained by refining
(\ref{pato2}) by means of the decomposition
\begin{equation}
\label{fulsi}
M(c)=\underset{P\in \P^*(c)}\oplus N(P).
\end{equation}
Here in general $\vert\P^*(c)\vert\neq 1$. In particular $L$ is not always uniserial.
This paper does not study whether $c(L)=\vert\P^*\vert$, and hence (\ref{fulsi}), hold for all $n$. As illustration
we refer the reader to Examples \ref{tatin1} and \ref{tatin2}.

Let us turn to the contents of the paper. Section \ref{pre} contains definitions and notation,
as well as basic facts about $I$ and $L$ to be used throughout the article. It also defines $N(P)$,
whose irreducibility is proven in Section \ref{nulam}.
Section \ref{cmp} associates to every $P\in\P(c)$ a suitable $P^*\in\P^*(c)$ satisfying $N(P^*)=N(P)$,
explicitly computes the common value $\nu_\ell([P:B])=\nu_\ell([P^*:B])$, and shows that
$M(c)$ is the direct sum of all distinct $N(P^*)$ as $P^*$ runs through $\P^*(c)$. These proofs are long.
To not interrupt the flow of the paper we created Appendix A to deal with the transfer of information from the lattice $\P$ of standard parabolic subgroups to the Steinberg lattice $I$ and hence to its modular reduction $L$,
and Appendix B to develop the auxiliary tools to find the exact value of $\nu_\ell([P:B])$.

Section \ref{cmp} also proves that the top
factor $M(0)=L(0)/L(1)$ is irreducible, where
$L(1)$ is the only maximal submodule of $L$, i.e. $rad(L)=L(1)$.
This result is dual to the aforementioned fact that $M(b)=L(b)=soc(L)$ is irreducible.
That this is not to be taken for granted is shown by Examples 5.4 and 5.5 of [4], where the reduction modulo $\ell$ of the Steinberg lattice
of $\Sp(4,q)$ is seen not to be irreducible modulo its radical.

We spend considerable effort -see Theorems \ref{A8} and \ref{radicalalfonsin}- demonstrating that the socle and radical series of $L$ simply agree with (\ref{pato2}), provided the positive integer
\begin{equation}
\label{dedd}
d=\nu_\ell\big(\frac{q^e-1}{q-1}\big)
\end{equation}
is equal to one. This is a pleasant state of affairs taking into account how differently these series are defined
and the fact that $L$ in general is not uniserial, even when $d=1$. We do not know if Theorems \ref{A8} and \ref{radicalalfonsin} hold when $d>1$.

Regarding $c(L)$, we know from above that $V\leq c(L)\leq |\P^*|$. Under strong hypotheses, such as in
Theorem \ref{A9}, all three of these numbers coincide. In general, we have
a recursive formula for $|\P^*|$ (Lemma \ref{longit}) and an explicit one for $V$ (Theorem \ref{A10}).
In most cases $V$ is a
polynomial in $\ell$ which depends on $d$ and the digits of
$\lfloor n/e \rfloor$ when written in base $\ell$.

Theorem \ref{thorx} finds a new generator for $I$ that is a common eigenvector when $U$ acts on~$I$.
While the statement of our result makes sense for all
finite groups of Lie type, it need not hold outside of type
$A$. Indeed, in Examples 5.4 and 5.5 of [4] we find  that the
top factor of the analogue of $L$ for $\mathrm{Sp}(4,q)$ is
completely reducible. If any common eigenvector for the action of
$U$ generated $I$, the top factor of $L$ would be irreducible,
against [4].

Theorem \ref{A12} computes the endomorphism ring of any term of the series (\ref{pato1}) or
(\ref{pato2}): it consists entirely of scalar operators.

As noted in [4], each term $I(c)$ of (\ref{pato1}) is free of rank $|U|$. It follows
that $I(c)$ is a Steinberg lattice. Let $T^c=I(c)/\ell I(c)$ stand for the reduction
of $I(c)$ modulo $\ell$. In this notation, $L=T^0$. The $KG$-modules $T^c$ are studied in Section \ref{scvsl}.
Surprisingly, they
are pairwise non-isomorphic for all $\P$-values $c$. Consequently, the $RG$-modules $I(c)$, for all $\P$-values $c$, are
non-isomorphic to each other. By a well-known theorem of Brauer and Nesbitt,
the non-isomorphic $KG$-modules $T^c$ must
have the same composition factors. We obtain a direct proof of this fact
by comparing the factors of the series
(\ref{pato1}) and (\ref{pato2}). We also find that the socle of $T^c$, $0<c<b$, is no longer irreducible, as it contains copies of $M(b)=L(b)$ and $M(0)$.
This contrasts with the cases of $L$ and $T^b\cong L^*$, both of which have an irreducible socle.
We wonder if there is a Steinberg lattice whose reduction modulo $\ell$ is completely reducible
for $\ell\mid [G:B]$.

Finally, we investigate when $\vert \P^*(c)\vert=1$ for all $\P$-values $c$, which, by above,
is a sufficient condition for (\ref{pato2}) to be a composition series of $L$.
The answer, Theorem \ref{A1},
depends on various cases and will not be stated here. Its case by case proof is given in Appendix C.
In any case, (\ref{pato2}) is a composition series of $L$
if $\lfloor n/e \rfloor\leq \ell$. The bottom $\ell+1$ non-zero factors of (\ref{pato2}) remain irreducible
if $\lfloor n/e \rfloor>\ell$ (Theorem \ref{A2}).

\section{Preliminaries}
\label{pre}


Let $e_1,...,e_n$ be the canonical basis of the column space
$F_q^n$. For $\sigma\in S_n$ we have the permutation matrix
$\widetilde{\sigma}\in G$ given by
$\widetilde{\sigma}e_i=e_{\sigma(i)}$. We abuse notation and
identify $\sigma$ with~$\widetilde{\sigma}$.

Let $\Pi$ be the set of all fundamental transpositions
$(1,2),...,(n-1,n)$. There is a natural bijection from the set of
all subsets of $\Pi$ onto $\P$, given by $J\mapsto P_J=\langle
B,J\rangle$.

To any $(i,j)$, with $1\leq i\neq j\leq n$, there corresponds the
root subgroup $X_{ij}$ of $G$ formed by all matrices
$t_{ij}(a)=I_n+aE^{ij}$, as $a$ runs through $F_q$.

Let $R^*$ stand for the unit group of $R$. To a group homomorphism $\l:U\to R^*$ we associate the set
$J(\l)\subseteq \Pi$ of all $(i,i+1)$ such that $\l$ is
non-trivial on $X_{i,i+1}$ and let
$$
P(\l)=P_{J(\l)}
$$
be the corresponding standard parabolic subgroup. Every $P\in \P$
arises in this way.

\subsection{$M(c)\neq 0$ for every $\P$-value $c$}

Fix a group homomorphism $\l:U\to R^*$
and let $P=\P(\l)\in\P$, $c=\nu_\ell([G:P])$. Associated to $\l$ we have the element $E_\l$ of $I$ defined by
\begin{equation}
\label{defegrande} E_\l=\underset{u\in U}\sum\l(u)u\u\in I.
\end{equation}
We see $U$ acts on $E_\l$ via $\l^{-1}$ and any $x\in I$ with this
property is a scalar multiple of $E_\l$.

Let $f$ be the bilinear form on $I$ defined in the Introduction.
As seen in Section 3 of [4]
\begin{equation}
\label{fatigue}
f(E_\l, u\u)=\l(u)[G:P],\quad u\in U.
\end{equation}
It follows that
\begin{equation}
\label{np5}
E_\l\in I(c)
\end{equation}
and
\begin{equation}
\label{np6}
E_\l\notin I(c+1).
\end{equation}
Let
$$
F_{\l}=E_\l+\ell I\in L.
$$
We see from above that
\begin{equation}
\label{np7}
F_\l\in L(c).
\end{equation}
It was asserted in Section 4 of [4] that
\begin{equation}
\label{estovale} F_\l\notin L(c+1).
\end{equation}
This does not follow automatically from (\ref{np6}), and
we pause to verify this crucial assertion. We need to show that $E_\l\notin I(c+1)+\ell I$.
Since $P(\l)=P(\l^{-1})$ we have
$E_{\l^{-1}}\in I(c)$, as above. Therefore, for all $x\in I(c+1)+\ell I$
$$
f(x,E_{\l^{-1}})\in \ell^{c+1}R.
$$
But (\ref{fatigue}) gives
$$
f(E_\l,E_{\l^{-1}})=|U|[G:P(\l)],
$$
where $\ell\nmid |U|$, so indeed $E_\l\notin I(c+1)+\ell I$, as claimed.
Combining (\ref{estovale}) and (\ref{np7}) we obtain
\begin{thm}
\label{gp3}
(Gow) The $KG$-module $M(c)=L(c)/L(c+1)\neq 0$ for every $\P$-value $c$.
\end{thm}

Since $F_\l\in L(c)$ but $F_\l\notin L(c+1)$ we see that
\begin{equation}
\label{cuw}
N(\l)=KG\cdot (F_\l+L(c+1)).
\end{equation}
is a non-zero cyclic submodule of $M(c)$. We will see shortly that $N(\l)$ is irreducible.

\subsection{$M(c)\neq 0$ implies that $c$ is a $\P$-value}

To derive the converse of Theorem \ref{gp3} we require two further tools. The first, taken from
from [6], was originally proven by Gelfand and Graev for complex~representations.

\begin{thm}
\label{gg} A non-zero $KG$-module has a one dimensional
$U$-invariant subspace.
\end{thm}

\begin{lem}
\label{todos} The natural group homomorphism $\l\mapsto
\overline{\l}$, where $\overline{\l}(u)=\l(u)+\ell R$, from the
group of all group homomorphisms $U\to R^*$ to the group of all
group homomorphisms $U\to K^*$, is an isomorphism.
\end{lem}

\noindent{\it Proof.} Since $U/U'$ is an elementary abelian
$p$-group and both $R^*$ and $K^*$ possess a non-trivial $p$-root
of unity, we see that the groups our map is connecting have the
same size, namely $|U/U'|$. It thus suffices to show that our map
is injective. For this purpose, suppose that $\overline{\l}$ is
trivial. We wish to show that $\l$ must be trivial. If not, then
$\l(u)=a\neq 1$ for some $u\in U$. As $\overline{\l}$ is trivial,
$x=a-1\in \ell R$. Thus $a=1+x$ is a $p$-root of unity with $x\neq
0$ in $\ell R$. Let $k\geq 1$ be the $\ell$-valuation of $x$. Then
the $\ell$-valuation of $x^p$ is $kp>k$. But
$$1=a^p=(1+x)^p=1+px+\cdots+px^{p-1}+x^p.$$
Subtracting 1 from each side yields $x^p=-px(1+c)$, where $c\in
\ell R$. Since $p$ and $1+c$ are units in $R$, we reach the
contradiction that the $\ell$-valuation of $x^p$ is $k$.

\begin{prop}
\label{tufo} Let $c\geq 0$ be arbitrary and let
$M$ be a $KG$-submodule of $L$ properly
containing $L(c+1)$. Then $M$ contains $F_\l$ for some
group homomorphism $\l:U\to R^*$ such that $k=\nu_\ell([G:P(\l)])\leq c$.
If actually $M\subseteq L(c)$ then $k=c$.
\end{prop}

\noindent{\it Proof.} By assumption $M/L(c+1)$ is a non-zero
$KG$-module. Then $M/L(c+1)$ has a one dimensional
$U$-invariant subspace, say $A/L(c+1)$, where $A$ is a
$KU$-submodule of $M$, by Theorem \ref{gg}. Since $\ell\nmid |U|$,
$A$ is completely reducible as a $KU$-module. Let $N$ be a
$KU$-complement to $L(c+1)$ in $A$. Then $N$ is a one
dimensional $KU$-submodule of $M$ not contained in $L(c+1)$.

Now $U$ acts on $N$ via a linear character, say $\mu:U\to K^*$.
From Lemma \ref{todos} we know that $\mu=\overline{\eta}$ for a
unique linear character $\eta:U\to R^*$. Let $\l=\eta^{-1}$. We easily see that $U$ acts
on $F_{\l}$ via $\mu$. Since $U$ acts on $L$ via the regular
representation, it follows that $N=K\cdot F_{\l}$. Thus $F_\l$
is in $M$ but not in $L(c+1)$.

Suppose, if possible, that $k>c$. Then $k\geq c+1$, so $L(k)\subseteq L(c+1)$. But
$F_\l\in L(k)$ by (\ref{np7}), so $F_\l\in L(c+1)$, a contradiction. This proves the first
assertion.

Assume next that $M\subseteq L(c)$. Suppose, if possible, that $k<c$. Then $k+1\leq c$ so
$F_\l \in M\subseteq L(c)\subseteq L(k+1)$,
against (\ref{estovale}). This completes the proof.

\medskip

Proposition \ref{tufo} applied to $M=L(c)$ yields

\begin{thm}
\label{xzc}
 (Gow) If the $KG$-module $M(c)=L(c)/L(c+1)\neq 0$
then $c$ is a $\P$-value.
\end{thm}

\subsection{Notation associated to parabolic subgroups}
\label{ghk}

Let $H$ be the diagonal subgroup of $G$. As $U$ is normalized by
$H$ we have an action of $H$ on the set of all group homomorphisms
$U\to R^*$. The orbits of this action are parametrized by $\P$.
Indeed, the $H$-orbit of $\l:U\to R^*$ is formed by all $\mu:U\to
R^*$ such that $P(\lambda)=P(\mu)$.

Fix $P\in\P$ for the remainder of this subsection and let $c=\nu_\ell([G:P])$.

Given a group homomorphism $\l:U\to R^*$ and $h\in H$ we see
from (\ref{defegrande}) that
$$
h\cdot E_\l=E_{\, ^{h}\!\l},\quad h\in H.
$$
If $\mu:U\to R^*$ is also group homomorphism and $P(\l)=P(\mu)$ we can then find $h$ in $H$ such that $^{h}\!\l=\mu$. We conclude that $RG\cdot E_\l=RG\cdot E_\mu$ whenever $P(\lambda)=P(\mu)$. We may
thus define $RG$-submodule $I'(P)$ of $I$ by
$$
I'(P)=RG\cdot E_\l
$$
for any $\l$ such that $P(\lambda)=P$. We also define the $KG$-submodule $L'(P)$ of $L$
by
$$
L'(P)=(I'(P)+\ell I)/\ell I=KG\cdot F_\l
$$
and the $KG$-submodule $N(P)$ of $M(c)$
$$
N(P)=(L'(P)+L(c+1))/L(c+1)=KG\cdot (F_\l+L(c+1))=N(\l)
$$
for any choice of $\l$ satisfying $P(\lambda)=P$. We further define
$$
I(P)=I(c),\quad L(P)=L(c),\quad M(P)=M(c).
$$
In this notation, we have the following result.
\begin{cor}
\label{nuevo} Eliminating repeated terms from
(\ref{pato2}) produces (all inclusions are proper):
\begin{equation}
\label{filt} 0\subset L(P_0)\subset\cdots\subset L(P_{V-1})=L,
\end{equation}
where $P_0,...,P_{V-1}\in\P$ are chosen so that
$\nu_\ell([G:P_0])>\cdots>\nu_\ell([G:P_{V-1}])$.
\end{cor}
We find it useful to have a notation to pass from one term of (\ref{filt})  to the next. Let
$$
L(P)^\sharp=L(c+1).
$$
Thus $L(P)^\sharp=0$ if $c=b$ and $L(P)^\sharp=L(Q)$
if $\nu_\ell([G:Q])$ is the first $\P$-value larger than $c$.

\section{Irreducibilty of $N(\l)$}
\label{nulam}

We quote the following result from [6].

\begin{thm}
\label{credu}
$M(c)$ is completely reducible and self-dual, while $L$ is multiplicity free.
\end{thm}

\begin{thm}
\label{yucatan}
Let $\l:U\to R^*$ be a group homomorphism
with $c=\nu_\ell([G:P(\l)])$. Then $N(\l)$, as defined in (\ref{cuw}), is an absolutely irreducible $KG$-sumodule of $M(c)$.
\end{thm}

\noindent{\it Proof.} We know from Theorem \ref{credu} that $M(c)$, and hence $N(\l)$, is completely reducible,
so it suffices to
show that the only $KG$-endomorphisms of $N(\l)$ are scalars.

Let $\mu:U\to K^*$ be the group homomorphism corresponding to $\l^{-1}:U\to R^*$ by the natural
projection $R^*\to K^*$. Since $U$ acts on $L$ via the regular representation and $\ell\nmid \vert U\vert$
we see that $\mu$ enters a given $KU$-section of $L$ at most once. But construction, if
$x=F_\l+L(c+1)$, then $u\cdot x=\mu(u)x$ for all $u\in U$. Moreover, $0\neq x\in N(\l)$ shown in Section \ref{pre}.
It follows that the subspace of $N(\l)$ where $U$ acts via $\mu$ is one dimensional and spanned by $x$.
Let $\alpha$ be an arbitrary $KG$-endomorphism of $N(\l)$. If $u\in U$ then
$$u\alpha(x)=\alpha(ux)=\alpha(\mu(u)x)=\mu(u)\alpha(x),$$
whence $\alpha(x)=a x$
for some $a\in K$ by above. But $N(\l)=KG x$ by construction, so
if $y\in N(\l)$ then $y=rx$ for some $r\in KG$, whence $\alpha(y)=ay$, as required.

\begin{cor}
\label{loriverde}
Let $c$ be a $\P$-value. Then every irreducible submodule of $M(c)$ must be of
the form $N(\l)$ for some $\l:U\to R^*$
satisfying $c=\nu_\ell([G:P(\l)])$.
\end{cor}

\noindent{\it Proof.} This follows from Proposition \ref{tufo}.

\begin{cor} All irreducible constituents of $L$ are absolutely irreducible.
\end{cor}

\noindent{\it Proof.} Using the series (\ref{pato2}), this follows from Theorem \ref{credu}
and Corollary \ref{loriverde}.

\section{Construction of $\P^*$ and first consequences}
\label{cmp}

A composition of $n$ is a sequence $(a_1,...,a_k)$ such that
$a_1,...,a_k$ are positive integers adding up to $n$. There is a
natural bijection from the set of all compositions of $n$ onto
$\P$, given by $(a_1,...,a_k)\mapsto P_{(a_1,...,a_k)}$, the block
upper triangular group with blocks of sizes $a_1,...,a_k$. By
abuse of notation we will identify each $P\in\P$ with its
corresponding composition.

A parabolic subgroup $Q=(b_1,...,b_l)$ is equivalent to $P$ if
$k=l$ and $(b_1,...,b_k)$ is a rearrangement of $(a_1,...,a_k)$.
Thus, the parabolic subgroups equivalent to $P$ can be obtained by
repeated application of single swaps of the form
$a_i\leftrightarrow a_{i+1}$.

\begin{thm}
\label{pringle}
If $Q\in\P$ is equivalent to a subgroup of $P\in\P$ then $I'(Q)\subseteq I'(P)$ and,
consequently, $L'(Q)\subseteq L'(P)$.
\end{thm}

\noindent{\it Proof.} This can be found in Appendix A.

\begin{cor}
\label{pringle2}
If $P$ and $Q$ are standard parabolic subgroups,
$\nu_\ell([G:P])=\nu_\ell([G:Q])$, and $Q$ is equivalent to a subgroup of $P$, then
$N(Q)=N(P)$.
\end{cor}

\noindent{\it Proof.} This follows from Theorems \ref{yucatan} and \ref{pringle}.

Let
$$
m=\mathit{max}\{i\geq 0\,|\, \ell^i\leq \lfloor n/e \rfloor\}.
$$
Given $1\leq a\leq n$ we write
$$
\Delta(a)=(y_{-1},y_0,\dots,y_m),
$$
where $0\leq y_{-1}<e$, $0\leq y_i<\ell$ for $1\leq i\leq m$, and
$a=y_{-1}+e(y_0+y_1\ell+\cdots+y_m\ell^m)$.
Thus $y_{-1}$ is the remainder of dividing $a$ by $e$ and
$(y_m\dots y_0)_\ell$ is the representation of $\lfloor
a/e\rfloor$ in base $\ell$. Given $P=(a_1,\dots,a_k)\in \P$ we let
$$
\Delta(P)=\Delta(a_1)+\cdots+\Delta(a_k).
$$
Thus $\Delta(P)=(z_{-1},z_0,\dots,z_m)$ is a sequence
non-negative integers satisfying
$$
z_{-1}+z_0e+z_1e\ell+\cdots+z_me\ell^m=n.
$$
We define
\begin{equation}
\label{demaria}
P^*=(\underbrace{1,\dots,1}_{z_{-1}},
\underbrace{e,\dots,e}_{z_0},\underbrace{e\ell,\dots,e\ell}_{z_1},\dots,
\underbrace{e\ell^m,\dots,e\ell^m}_{z_m})=[z_{-1},z_0,\dots,z_m].
\end{equation}
Let $\P^*$ be set of all standard parabolic subgroups of the
form (\ref{demaria}). They correspond to partitions of $n$ where each part is
either 1 or of the form $e\ell^i$ for some $0\leq i\leq m$.

Recall the definition (\ref{dedd}) of $d$. We define the sequence $s_0,s_1,...$ of positive integers by
$$
s_0=d, s_1=\ell d+1, s_2=\ell^2d+\ell+1, s_3=\ell^3d+\ell^2+\ell+1,\dots.
$$

\begin{thm}
\label{muyjota} Let $P=(a_1,\dots,a_k)\in \P$.  Then $P^*$, as defined in (\ref{demaria}), is
equivalent to a parabolic subgroup contained in $P$. Moreover,
$$
\nu_\ell([P:B])=s_0z_0+\cdots+s_mz_m=\nu_\ell([P^*:B])\;\text{ and }\; N(P)=N(P^*).
$$
\end{thm}

\noindent{\it Proof.} The very construction of $P^*$ yields the
first assertion. It is known and easy to see that
$$
[P:B]=\underset{1\leq i\leq k}\Pi\;\;\underset{1\leq
j\leq a_i}\Pi (q^j-1)/(q-1).
$$
Let us write $\Delta(a_i)=(y^{i}_{-1},y^{i}_0,\dots,y^{i}_m)$. Then by Lemma \ref{r4} of Appendix B we have
$$
\nu_\ell([P:B])= \underset{1\leq i\leq k}\sum
\nu_\ell\big(\underset{1\leq j\leq a_i}\Pi (q^j-1)/(q-1)\big)
=\underset{1\leq i\leq k}\sum (y^{i}_m
s_m+\cdots+y^{i}_0s_0)=z_ms_m+\cdots+ z_0s_0.
$$
As for $P^*$, the same argument (but using Lemma \ref{v8} of Appendix B instead) yields
$$
\nu_\ell([P^*:B])=\underset{0\leq i\leq m}\sum z_i\times
\nu_\ell\big(\underset{1\leq j\leq e\ell^i}\Pi (q^j-1)/(q-1)\big)
=z_0s_0+\cdots+z_m s_m.
$$
This proves the second assertion. We may now derive the third from Corollary \ref{pringle2}.

\begin{cor}
\label{saporeve}
Let $\lfloor n/e \rfloor=(x_m\dots x_0)_\ell$. Then $b=\nu_\ell([G:B])=s_0x_0+\cdots+s_mx_m$.
\end{cor}

\begin{note} Let $P\in\P$. Then $P^*$ is the only member of $\P^*$ that is
equivalent to a standard parabolic subgroup contained in $P$ and
satisfies $\nu_\ell([P:B])=\nu_\ell([P^*:B])$.
\end{note}

\begin{thm}
\label{niatro}
Let $c$ be a $\P$-value. Then $M(c)$ has the following decomposition into non-isomorphic irreducible $KG$-modules:
$$
M(c)=\oplus N(P),
$$
where the sum runs through all different $N(P)$ produced by the $P\in \P^*(c)$.
\end{thm}

\noindent{\it Proof.} This follows from Theorems \ref{credu}, \ref{yucatan}, \ref{muyjota} and Corollary \ref{loriverde}.

\begin{cor}
\label{icrit} $M(c)$ is irreducible if and
only if $c$ is a $\P$-value and $N(P)=N(Q)$ for all $P,Q\in\P^*(c)$. In particular,
if $|\P^*(c)|=1$ then $M(c)$ is irreducible, and if $|\P^*(c)|=1$ for $\P$-values $c$ then
(\ref{pato2}) is a composition series of $L$.
\end{cor}

\begin{thm}
\label{thorx}
Let $\l:U\to R^*$ be a group homomorphism such that $P(\l)=G$, i.e. $\l$ is non-trivial in every fundamental
root subgroup. Then $E_\l$ and $F_\l$, as defined in Section~\ref{pre}, respectively generate $I$ and $L$, i.e.
$I=RG\cdot E_\l$ and $L=KG\cdot F_\l$.
\end{thm}

\noindent{\it Proof.} Let $M$ the $KG$-submodule of $L$ generated by the $F_\mu$ as $\mu$ runs
through all group homomorphisms $U\to R^*$. By hypothesis $P(\mu)\subseteq P(\l)$ for every $\mu$, so Theorem \ref{pringle} gives $KG\cdot F_\l=M$. On the other hand, Theorem \ref{gg} and Lemma \ref{todos}
show that $M=L$. Hence $KG\cdot F_\l=L$. Since the image of $E_\l$ generates $L=I/\ell I$ and $I$ is a finitely generated $R$-module, it follows from Nakayama's lemma that $E_\l$ generates $I$.

\begin{note} If $P(\l)=G$ then Lemma \ref{ayuda} of Appendix A gives
$E_\l=sg(\sigma_0)\underset{u\in U}\sum \l(u)u\sigma_0 \widehat{B}$, where $\sigma_0\in S_n$
is defined by $\sigma_0=(1,n)(2,n-1)(3,n-2)\cdots$.
\end{note}

\begin{thm}
\label{A3} The top factor $M(0)=L(0)/L(1)$ of (\ref{pato2}) is always irreducible.
\end{thm}

\noindent{\it Proof.} $M(0)=N(G)$ by Theorem \ref{thorx}, so Theorem \ref{yucatan}
applies.

\begin{thm}
\label{maxl1}
All proper submodules of $L$ are contained in $L(1)$, i.e. $rad(L)=L(1)$.
\end{thm}

\noindent{\it Proof.} By above $L(1)$ is a maximal submodule of $L$.
Suppose, if possible, that $M$ is a proper submodule of $L$ different from $L(1)$. Therefore $L(1)+M=L$, so
by the second isomorphism theorem $L(1)/L(1)\cap M\cong (L(1)+M)/M=L/M$. Let $\l$ be a linear character of $U$ that is non-trivial in every fundamental root subgroup. As $L/L(1)\neq 0$ by Theorem~\ref{gp3} and
$L/M\neq 0$ by assumption, Theorem \ref{thorx} implies that
$\l$ enters both $L/L(1)$
and $L/M$. Hence $\l$ enters both factors of the series $L\supset L(1)\supset L(1)\cap M$. This
contradicts the fact that $U$ acts on $L$ via the regular representation with $\ell\nmid |U|$.

\section{Endomorphism rings of $I(c)$ and $L(c)$}

\begin{lem} Every $L(P)$ is equal to the sum of the submodules $L'(Q)$ inside it.
\label{insidet}
\end{lem}

\noindent{\it Proof.} This is certainly true for the irreducible
module $L(B)$. Suppose that $\nu_\ell([P:B])>0$ and the statement is true
$L(P)^\sharp$. Since $L(P)$ equals $L(P)^\sharp$ plus the sum of
certain $L'(Q)$ inside $L(P)$ by Theorem \ref{niatro}, the result
follows by induction.

\begin{lem}
\label{mickey} Let $P\in\P$. Then the only endomorphisms of $L(P)$ are scalars.
\end{lem}

\noindent{\it Proof.} Let $f$ be an endomorphism of $L(P)$. Then $f$ is determined by
its values on the submodules $L'(Q)$ inside $L(P)$ by Lemma \ref{insidet}.
Arguing as in the proof of Theorem \ref{yucatan} we see that if $\l:U\to K^*$ is any group homomorphism
then $f(F_\l)=a(\l)F_\l$ for some $a(\l)\in K$. Let $\l_0:U\to K^*$ be the trivial group homomorphism.
Then $F_{\l_0}$ belongs to all $KG\cdot F_{\l}$ by Theorem
\ref{pringle}. It follows that all $a(\l)$ are equal to $a(\l_0)$, as required. (If
$P=G$ we may use Theorem \ref{thorx} to simplify the above argument.)

\begin{thm}
\label{A12}
Each endomorphism of a term of the series (\ref{pato1}) or (\ref{pato2}) is a scalar.
\end{thm}

\noindent{\it Proof.} Let $F$ be the field of fractions of $R$. Then the
Steinberg module $F\otimes_R I(c)$ over $F$ is absolutely
irreducible, so its only endomorphisms are scalars. It follows
that the only endomorphisms of $I(c)$ are scalars. The case of $L(c)$ is given in Lemma \ref{mickey}.

\begin{note}
\label{fillol} Every $I(c)$ is generated as an $RG$-module by
elements of form $\ell^i E_\l$.

Indeed, this is true for $I(0)$ by Theorem
\ref{thorx}. Suppose $c>0$ and the result is true for $I(c-1)$. Let
$N$ be the sum of all submodules $RG\cdot E_\l$ inside $I(c)$ and
consider the $KG$-module $M=I(c)/(N+\ell I(c-1))$. We wish to show
that $M=0$. Consider the natural epimorphism $I(c)\to L(c)$. Its
kernel is $I(c)\cap \ell I=\ell I(c-1)$. Thus $I(c)/\ell
I(c-1)\cong L(c)$. Under this isomorphism $N+\ell I(c-1)/\ell
I(c-1)$ corresponds to the submodule of $L(c)$ generated by all
$F_\l$ inside $L(c)$, namely $L(c)$,  by Lemma \ref{insidet}.
Thus
$$M/(N+\ell I(c-1))\cong (I(c)/\ell I(c-1))/(N+\ell I(c-1)/\ell I(c))\cong L(c)/L(c)=0.$$
Therefore $M=N+\ell I(c-1)$,
and the result follows by induction.
\end{note}

\section{Socle and Radical series of $L$}
\label{noahuant}

For $P\in\P$ we set
$$
\vartheta(P)=\nu_\ell([G:P]),\quad \p(P)=\nu_\ell([P:B]).
$$
Using $[G:B]=[G:P][P:B]$ and that $V$ is the total number of $\P$-values we find that
$$
|\{\vartheta(P)\,:\, P\in\P\}|=V=|\{\p(P)\,:\, P\in\P\}|.
$$
Recall that $e=\mathit{min}\{\,i\geq 2\,|\, \ell\text{ divides }\frac{q^i-1}{q-1}\}$ and set $m=\mathit{max}\{i\geq 0\,|\, \ell^i\leq \lfloor n/e \rfloor\}$.
Let
$$
\lfloor
n/e \rfloor=(x_m\dots x_0)_\ell=x_m\ell^m+\cdots+x_1\ell+x_0.
$$
Given nonnegative integers $z_0,z_1,\dots,z_m$ satisfying
$e(z_0+\cdots+z_m\ell^m)\leq n$ we set
$z_{-1}=n-e(z_0+z_1\ell+\cdots+z_m\ell^m)$ and reduce the notation
$[z_{-1},z_0,z_1,\dots,z_m]$ of (\ref{demaria}) to $[z_0,z_1,\dots,z_m]$.
Thus,
$$
[z_0,z_1,\dots,z_m]=(\underbrace{1,\dots,1}_{z_{-1}},
\underbrace{e,\dots,e}_{z_0},\underbrace{e\ell,\dots,e\ell}_{z_1},\dots,
\underbrace{e\ell^m,\dots,e\ell^m}_{z_m})\in\P^*.
$$
Recall also that $d=\nu_\ell\big(\frac{q^e-1}{q-1}\big)$ and $b=\vartheta(B)$, whose
exact value is given in Corollary \ref{saporeve}. We will also appeal to
the notation introduced in Section \ref{pre}.

\begin{thm}
\label{A2}

(a) If $\lfloor n/e\rfloor\leq\ell$
then $L$ is uniserial and its only composition series
is~(\ref{pato2}).

(b) If $\lfloor n/e\rfloor>\ell$ then the first $\ell+1$ terms of the socle
series of~$L$, together with 0, are $0\subset
L(P_0)\subset\cdots\subset L(P_\ell)$, in the notation of Corollary \ref{nuevo}.
This is in fact a composition series of $L(P_\ell)$. In particular, $L(P_\ell)$ is
uniserial of length $\ell+1$ and the first $\ell+1$ factors of
(\ref{filt}) starting from the bottom are irreducible.
\end{thm}

\noindent{\it Proof.} (a) Note that if $\lfloor
n/e\rfloor<\ell$ then $\P^*=\{[i]\,|\, 0\leq i\leq \lfloor
n/e\rfloor\}$, while if
 $\lfloor n/e\rfloor=\ell$ then
 $\P^*=\{[i,0]\,|\, 0\leq i\leq \ell\}\cup\{[0,1]\}$. In both
 cases $\P^*$ is ordered by inclusion, which explains why $L$ is
 uniserial.

Indeed, let us agree that the socle series of $L$ starts at 0. Let
$P\in\P^*$. Suppose that $L(P)^\sharp$ is equal to
 a term of the socle series of $L$ and let $S$ be the next term of
 this series. We wish to show that $S=L(P)$ with $S/L(P)^\sharp$
 irreducible.

We have $L(P)\subseteq S$ by Theorem \ref{credu}. Let $M$ be
 a submodule of $L$ properly containing $L(P)^\sharp$ with
 $M/L(P)^\sharp$ irreducible. We know from Proposition \ref{tufo} that $M$ contains $L'(Q)$ for
some $Q\in\P^*$ satisfying $\vartheta(Q)\leq \vartheta(P)$. As
$\P^*$ is ordered by inclusion, $\vartheta(Q)\leq \vartheta(P)$
implies that $P$ is contained in $Q$. This implies $L'(P)\subseteq
L'(Q)$ by Theorem \ref{pringle}. Thus $M/L(P)^\sharp$ contains
$N(P)$, so $M/L(P)^\sharp=N(P)$ by the
irreducibility of $M/L(P)^\sharp$. As $M$ was arbitrary, $S/L(P)^\sharp$ itself is irreducible and equal to
$N(P)$. In particular, $S\subseteq L(P)$.

(b) Let ${\cal R}=\{[i,0,\dots,0]\,|\, 0\leq i\leq \ell\}$. It is
easy to see that if $P\in {\cal R}$, $Q\in\P^*$ and
$\vartheta(Q)\leq \vartheta(P)$ then $P$ is equivalent to a
parabolic subgroup contained in $Q$. We may now repeat the above
proof with every $P\in{\cal R}$.

\begin{cor} (a) If $\lfloor n/e\rfloor\leq \ell$ then $L(P)=L'(P)$
is cyclic for all $P\in\P$.

(b) $L(P)=L'(P)$ is cyclic for all $P=[i,0,\dots,0]$, $0\leq i\leq \ell$.
\end{cor}

\noindent{\it Proof.} This follows from Theorem \ref{A2}, since in a
uniserial module every term of the socle series is generated by
any element not belonging to the previous term.

\begin{thm}
\label{A10} Let
$$
A=x_m(\ell^m+\cdots+\ell+1)+x_{m-1}(\ell^{m-1}+\cdots+\ell+1)+\cdots+x_1(\ell+1)+x_0+1,
$$
$$
Z=x_m(d\ell^m+\cdots+\ell+1)+x_{m-1}(d\ell^{m-1}+\cdots+\ell+1)+\cdots+x_1(d\ell+1)+x_0d+1,
$$
$$
C=x_m(\ell^{m-1}+\cdots+\ell+1)+x_{m-1}(\ell^{m-2}+\cdots+\ell+1)+\cdots+x_2(\ell+1)+x_1+1,
$$
noting that $C=1$ if $m=0$. Then

(a) $V=C+X$, where $X$ is the amount of values
$\phi(Q)$ satisfying $0\leq \phi(Q)<d\lfloor
n/e\rfloor$. Moreover, $X\geq \lfloor
n/e\rfloor$, so $A\leq V\leq Z$.

(b) Suppose $\lfloor n/e\rfloor\geq d\ell$. Then $V=Z-d^2\ell+Y$,
where $Y$ is the total amount of values $\phi(Q)$ satisfying
$0\leq \phi(Q)<d^2\ell$. Moreover, $Y\geq \ell d(d+1)/2$, so
$Z-\ell d(d-1)/2\leq V\leq Z$. In fact, if $d\leq \ell$ then
$Y=\ell d(d+1)/2$, so $V=Z-\ell d(d-1)/2$, that is
$$
V=x_m(d\ell^m+\cdots+\ell+1)+x_{m-1}(d\ell^{m-1}+\cdots+\ell+1)+\cdots+x_1(d\ell+1)+x_0d+1-\ell d(d-1)/2.
$$
\end{thm}

\noindent{\it Proof.} We may replace $\P$ by $\P^*$ in the statement of the theorem in view of Theorem \ref{muyjota}.

We will create a sequence of parabolic subgroups in $\P^*$
starting at $G^*=[x_0,x_1,...,x_m]$ and ending at
$[\lfloor n/e\rfloor,0,...,0]$. Our sequence will satisfy the following
properties: if $P$ is a term of the sequence and $P'$ is the next
term then $P'\subset P$ and $\phi(P')=\phi(P)-1$. The
number of terms of the sequence will be $C$. We will use Theorem \ref{muyjota} throughout.

The construction is as follows. Let $P\in\P^*$ and suppose $P$ is
not of the form $[a,0,...,0]$. Then
$P=[y_0,...,y_i,y_{i+1},0,...,0]$, where $0\leq i<m$ and
$y_{i+1}\neq 0$. We then define
$P'=[y_0,...,y_i+\ell,y_{i+1}-1,0,...,0]$. Starting at $G^*$ and
repeating this process $x_m$ times we reach
$[x_0,x_1,...,x_{m-1}+x_m\ell,0]$. Repeating now the process
$x_{m-1}+x_m\ell$ times we reach
$[x_0,x_1,...,x_{m-2}+x_{m-1}\ell+x_m\ell^2,0,0]$, and so on. All in all, our process
produces $C$ consecutive values, from $Z-1=\phi(G)$ to $d\lfloor n/e\rfloor=\phi([\lfloor n/e\rfloor,0,...,0])$.
This explains (a).

Suppose now $\lfloor n/e\rfloor\geq d\ell$. Given $P=[a,0,...,0]$,
where $d\ell< a\leq \lfloor n/e\rfloor$, define
$P^0=[a-1,0,...,0]$. Through a second process we can attain all
$d$ numbers from $da$, excluded, down to $d(a-1)$, included, as
values $\phi(Q)$. Given such $P=[a,0,...,0]$ define
$P^1=[a-(d\ell+1),d,0,...,0]$. Then $\phi(P^1)=\phi(P)$, and we
can now apply the first process $d$ times to $P^1$ until $P^0$ is
reached. Combining this with the above process, all numbers from
$\phi(G)$ down to $\phi([d\ell,0,...,0])$ are attained as values
$\phi(Q)$. This creates
$$
 C+d(\lfloor n/e\rfloor-d\ell)=Z-d^2\ell
$$
consecutive values $\phi(Q)$. This explains the first sentence of (b). We next show
that $Y\geq \ell d(d+1)/2$. Indeed, if $0\leq j\leq d-1$ and  $j\ell\leq
a<d\ell$  then $0\leq ad+j<d^2\ell$ is attained at
$Q=[a-j\ell,j,0,...,0]\in \P^*$. Thus $Y\geq \ell d+\ell(d-1)+\cdots+\ell=\ell d(d+1)/2$,
confirming the second sentence of (b). Next we
show $Y=d(d+1)/2$ provided $d\leq\ell$. We wish to know when a number
$0\leq h<d^2\ell$ is of the form $\phi(Q)$.
Now $d\ell^2+\ell+1\geq d^2\ell+\ell+1>d^2\ell$, so any such $Q$
will have to have the form $Q=[x,y,0,...,0]$.
Dividing $h$ by $d$, we may write $h=ad+j$, where
$0\leq a<d\ell$ and $0\leq j\leq d-1$. We look for $x,y$
such that $ad+j=\phi(Q)=dx+y(d\ell+1)$. Congruence modulo
$d$ reveals that $y\equiv j\mod d$. But if $y\geq d$ then
$y(d\ell+1)>d^2\ell$. Thus $y=j$. This implies $a=x+j\ell$, so
$a\geq j\ell$. The only attained values are the ones described above,
which completes the proof of (b).

\begin{note}
If $d\leq \ell$ and $\lfloor n/e\rfloor\leq d\ell$ the value of $V=c(L)$ is given in Theorem \ref{A9}.
This completes the determination of $V$ in all cases where $d\leq \ell$. Observe that if
$\lfloor n/e\rfloor=d\ell$ and $d\leq \ell$ then Theorems \ref{A9} and \ref{A10} compute $V$ in
different ways, but the answers agree. Indeed, if $d=\ell$ both give $V=\ell^3/2+ \ell^2/2+\ell+2$,
while if $d<\ell$ the common value is $\ell d(d+1)/2+d+1$.

The proof of Theorem \ref{A10} shows that at least $C-1$ {\em consecutive} top factors of (\ref{pato2})
are not zero. When $\lfloor n/e\rfloor\geq d\ell$ at least $Z-1-d^2\ell$ of them are non-zero. If
$d=1$ then $A=v=Z=b+1$ and {\em all} factors $L(c)/L(c+1)$, $0\leq c\leq b$, are non-zero.
\end{note}

\begin{thm}
\label{A8}
(a) If $d=1$ then (\ref{pato2}) is the socle series of $L$. (The 0 module and all prior repeated terms of
(\ref{pato2}) must be removed when interpreting this statement)

(b) Let $P\in\P$. Then $soc(L/L(P)^\sharp)=L(P)/L(P)^\sharp$,
except for the possibility that $(L'(Q)+L(P)^\sharp)/L(P)^\sharp$ be also irreducible,
where  $Q=[a,0,...,0]$, $\ell\leq a-1< \lfloor n/e\rfloor$ and $d(a-1)<\phi (P)<da$.
In particular, $soc(L/L(P)^\sharp)=L(P)/L(P)^\sharp$ if $\phi(P)\geq d\lfloor n/e\rfloor$.
 \end{thm}

\noindent{\it Proof.} Since $L/L(P)^\sharp$ is completely reducible, we always have
$L(P)/L(P)^\sharp\subseteq soc(L/L(P)^\sharp)$. We also know that $L(B)=soc(L)$, so equality
holds for $P=B$.

Suppose $\phi(P)>0$. By Proposition \ref{tufo} an arbitrary irreducible submodule of $L/L(P)^\sharp$ must have the form
$M=(L'(Q)+L(P)^\sharp)/L(P)^\sharp$, where $Q\in\P^*$ and $\phi(Q)\geq\phi(P)$.
Choose $Q$ so that $\phi(Q)$ is as small as possible. If
$\phi(Q)=\phi(P)$ then $M\subseteq L(P)/L(P)^\sharp$.

Suppose, if possible, that $\phi(Q)>\phi(P)$. If $d=1$ the
proof of Theorem \ref{A10} shows that $Q$ contains a parabolic
subgroup $Q'$ such that $\phi(Q')=\phi(Q)-1\geq\phi(P)$. By Theorem
\ref{pringle} we have $L'(Q')\subseteq L(Q)$, so the minimality of $Q$ is violated. Therefore, if $d=1$ we must have $\phi(Q)=\phi(P)$,
whence $M\subseteq L(P)/L(P)^\sharp$ and a fortiori $L(P)/L(P)^\sharp=soc(L/L(P)^\sharp)$. Induction then gives (a).
If $d$ is now arbitrary the proof of Theorem \ref{A10} yields the same contradiction as long as $Q$ is
not of the form $[a,0,...,0]$. Thus the case $\phi(Q)>\phi(P)$ can only occur when $Q=[a,0,...,0]$.
Since $Q\in \P$ we must have $a\leq \lfloor
n/e\rfloor$. As $\p(P)<\p(Q)$ we also have $0<a$ and $\p(P)<da\leq  d\lfloor
n/e\rfloor$. Now $[a-1,...,0]\subset Q$, so using Theorem
\ref{pringle} once more yields $\phi(P)>(a-1)d$. Since below $d\ell$
all values taken by $\phi$ decrease by $d$, it is also clear that
$a-1$ must be at least $\ell$. This proves (b).
\begin{prop}
\label{uniprop} (cf. Lemma \ref{insidet}) Suppose $d=1$. Then for every $\P$-value $c$, the $KG$-module $L(c)$ is the sum of all $L'(P)$ with
$P\in \P^*(c)$.
\end{prop}

\noindent{\it Proof.}  The result is true for the irreducible module $L(b)=L'(B)$. Suppose $c$ is a $\P$-value
smaller than $b$ and the result is true for the first $\P$-value $a$ larger than $c$.

We know from Theorem \ref{niatro} that $L(c)$ equals $L(a)$ plus the sum of submodules
$L'(P)$ such that $P\in\P^*(c)$. By inductive hypothesis, $L(a)$ is the sum of
all $L'(Q)$ such that $Q\in\P^*(a)$. Let $Q\in P^*(a)$. By Theorem \ref{pringle},
it suffices to find $P$ such that $P\in\P^*(c)$ and $Q$ is equivalent to a parabolic
subgroup contained in $P$. Let us write $Q=[y_0,...,y_m]$. If for any $i<m$ we have $y_i>\ell$
we can let $P$ be obtained from $Q$ by replacing $y_i$ by $y_i-\ell$ and $y_{i+1}$ by
$y_{i+1}+1$. We may therefore assume in what follows that $y_i<\ell$ for all $i<m$.
Let $y_{-1}=n-e(y_0+y_1\ell+\cdots+y_m\ell^m)$ and $x_{-1}=n-e \lfloor n/e \rfloor$.
We wish to show that $y_{-1}\geq e$, in which case
$P=[y_0+1,y_1,...,y_m]\in\P$. Using the hypothesis $d=1$ at this single point in the entire proof
ensures that $P$ satisfies our requirements.

We proceed to show that $y_{-1}\geq e$. Recalling that $\lfloor n/e \rfloor=(x_m\dots x_0)_\ell$, we first note that $y_m\leq x_m$. Indeed, if $y_m>x_m$ then $\ell^m y_m\geq \ell^m x_m+\ell^m$.
Using that all $x_j\leq \ell-1$ we easily see that $\ell^m>x_0+\cdots+x_{m-1}\ell^{m-1}$. Combining
these inequalities yields $\ell^m y_m>x_0+\cdots+x_{m-1}\ell^{m-1}+x_m\ell^m=\lfloor n/e \rfloor$,
so $\ell^m y_m\geq \lfloor n/e \rfloor+1$, whence $e\ell^m y_m\geq e(\lfloor n/e \rfloor+1)>n$,
contradicting the fact that $Q\in\P$.

Since $a>c\geq 0$, Theorem \ref{muyjota} eliminates the possibility that
$(y_m,\dots, y_0)=(x_m,\dots,x_0)$. Scan these sequences from left to right and let $i$ be the first index
satisfying $x_i\neq y_i$. Now argue as above, using $y_m=x_m,...,y_{i+1}=x_{i+1}$, to see that $y_i>x_i$
is impossible, so $y_i<x_i$. Suppose by way of contradiction that $y_{-1}< e$.
Now
$$0=(x_{-1}+e(x_0+\cdots+x_{i-1}\ell^{i-1}+x_i\ell^i))-(y_{-1}+e(y_0+\cdots+y_{i-1}\ell^{i-1}+y_i\ell^i)).$$
The largest possible value for the second summand is
$$
(e-1)+e(\ell-1+\cdots+(\ell-1)\ell^{i-1}+(x_i-1)\ell^i),
$$
namely
$$
(e-1)+e(\ell-1)(1+\cdots+\ell^{i-1})+e(x_i-1)\ell^i=
(e-1)+e(\ell^i-1)+e(x_i-1)\ell^i=ex_i\ell^i-1,
$$
while the smallest possible value for the first summand is $ex_i\ell$. This absurdity shows that $y_{-1}\geq e$,
thereby completing the proof.

\begin{thm}
\label{radicalalfonsin} If $d=1$ then (\ref{pato2}) (all repetitions removed) is the radical series of $L$.
\end{thm}

\noindent{\it Proof.} By convention $rad^0(L)=L$. Suppose $L(P)$ is a term of the radical series of $L$. We wish to
show that $rad(L(P))$ is $L(P)^\sharp$ (this will give a slightly different proof of Theorem \ref{maxl1}). Since $L(P)/L(P)^\sharp$ is completely reducible, it follows that $L(P)^\sharp$ contains $rad(L(P))$.
Suppose by way of contradiction that the inclusion is proper. Then the non-zero $KG$-module
$M=L(P)^\sharp/rad(L(P))$ must have a linear character $\l$ of $U$. By (\ref{estovale}) any such $\l$ must satisfy $\t(P(\l))>\t(P)$.
Now $L(P)/rad(L(P))$ is completely reducible, so its submodule $M$ is also a factor. Thus $M$ is a factor of $L(P)$.
By Proposition \ref{uniprop}, any non-zero image of $L(P)$ must necessarily contain a linear character $\mu$ of $U$ such that $\t(P(\mu))=\t(P)$. As remarked above, $M$
does not contain any such $\mu$, a contradiction.

\section{Comparing various Steinberg lattices}
\label{scvsl}

\begin{thm}
\label{A13} (a) $I(c)/I(c+1)\cong (L/L(c+1))^*$ ($*=\text{dual}$) as $KG$-modules, $c\geq 0$.

(b) For any $c\geq 0$ the composition factors of the $KG$-module $I(c)/I(c+1)$ are the composition factors of
$L(0)/L(1),...,L(c)/L(c+1)$ taken together.

(c) $I(0)/I(1)\cong L(0)/L(1)$.

(d) $I(c)/I(c+1)\cong L^*$ as $KG$-modules for all $c\geq b$.
\end{thm}

\noindent{\it Proof.} (a) Identify $K$ with $\ell^c R/\ell^{c+1}R$ and
consider the $RG$-homomorphism from $I(c)$ to $(L/L(c+1))^*$ given by $x\mapsto \varphi_x$, where $\varphi_x((y+\ell I)+L(c+1))=f(x,y)+\ell^{c+1}R$
for all $x\in I(c)$, $y\in I$, and $f$ is the bilinear form previously defined on $I$. Using results
from Section 4 of [4] we see that our map has kernel $I(c+1)$ and that $I(c)/I(c+1)$
and $L/L(c+1)$ have the same dimension, as required.

(b) The composition factors of $(L/L(c+1))^*$ are dual to those
of $L/L(c+1)$ (in reversed order). The composition factors of $L/L(c+1)$ are those
of $L(0)/L(1),...,L(c)/L(c+1)$, taken together, and these are all self-dual, so the result
follows from (a).

Alternatively, for $0\leq i\leq c$ there is a natural $RG$-epimorphism from $I(c+1)+\ell^i I(c-i)$ to $L(c-i)/L(c-i+1)$,
with kernel $I(c+1)+\ell^{i+1}I(c-(i+1))$ if $i\leq c-1$ and $I(c+1)$ if $i=c$.

(c) As $L/L(1)$ is self-dual, this follows from (a). Alternatively,
the natural epimorphism $I\to L/L(1)$ has kernel $\ell I+I(1)=I(1)$.

(d) Since $L(c+1)=0$ for all $c\geq b$, this follows from (a).

\begin{thm}
\label{B52}
Let $0\leq c<h$. Suppose there is a $\P$-value $a$ such that $c<a\leq h$. (In particular,
this applies when $h$ is a $\P$-value and when $c<h$ are both $\P$-values.)
Consider the Steinberg lattices $I(c)$ and $I(h)$
and let $T^c=I(c)/\ell I(c)$ and $T^h=I(h)/\ell I(h)$ be their respective reductions modulo $\ell$.
Then the $KG$-modules $T^c$ and $T^h$ are not isomorphic. Consequently, the $RG$-modules
$I(c)$ and $I(h)$ are not isomorphic.

On the other hand, $I(b)\cong I(b+1)\cong I(b+2)\cong\cdots$, so $T^b\cong T^{b+1}\cong T^{b+2}\cong \cdots$.
\end{thm}

\noindent{\it Proof.} The maps $v\mapsto \ell v\mapsto \ell^2v\mapsto\cdots$ yield
isomorphisms $I(b)\cong I(b+1)\cong I(b+2)\cong\cdots$, thereby justifying the last assertion.
Thanks to it, we may assume that $h\leq b$. We choose the $\P$-value $a$ to be as large as possible subject
to $a\leq h$.

We have
$$
I(c+1)/\ell I(c)=I(c+1)/I(c+1)\cap \ell I\cong (I(c+1)+\ell I)/\ell I=L(c+1),
$$
while by Theorem \ref{A13}
$$(I(c)/\ell I(c))/(I(c+1)/\ell I(c))\cong I(c)/I(c+1)\cong (L/L(c+1))^*.$$
Thus $T^c$ has a submodule isomorphic to $L(c+1)$ and the corresponding factor is
isomorphic to~$(L/L(c+1))^*$. The analogous result is valid
for $T^h$. Suppose there is an isomorphism from $T^c$ into $T^h$.
Now $L(c+1)$, and hence $T^c$, has a submodule
isomorphic to $L(b)$. Likewise, $T^h$ has a submodule isomorphic
to $L(b)$, unless $h=b$, in which case we must omit this part of
the proof and proceed to the next paragraph. Now the
$\ell$-modular reduction of any Steinberg lattice is multiplicity
free. Indeed, this just depends on the following facts: $U$ acts
on it via the regular representation; $\ell\nmid |U|$; any
non-zero $KG$-module must have a common eigenvector for $U$. Since $M(b)=L(b)$ is completely reducible, it follows that the supposed isomorphism must map the one copy of $L(b)$ inside $T^c$
into the one copy of $L(b)$ inside $T^h$. This induces an
isomorphism between the corresponding quotients. This process can
be continued.

Eventually, we get an isomorphism between a module $X$ with a
submodule isomorphic to $L(c+1)/L(h+1)$ with factor isomorphic to
$(L/L(c+1))^*$, and a module $Y$ isomorphic to $(L/L(h+1))^*$. If
$h$ is a $\P$-value then $a=h$, whereas if $h$ is not a $\P$-value
then $L(a+1)=\cdots=L(h+1)$. In any case, we may replace $h$ by
$a$ in the previous sentence.

Now $X$ has a submodule isomorphic to $M(a)$. But $Y$ does not have such a submodule.
For if it did, the dual of $Y$, namely $L/L(a+1)$, would have a factor isomorphic to the self-dual module
$M(a)$. Then $L$ would have the completely reducible module $M(a)$ as image. But $M(a)\neq 0$, since $a$ is a $\P$-value,
and $L$ has only one non-zero completely reducible image, up to isomorphism, namely the irreducible module $M(0)=L/L(1)$, as the radical $L(1)$ of $L$ is maximal. It would follow that $M(0)\cong M(a)$,
which is impossible since $a>0$ and $L$ is multiplicity free.

\begin{note} As mentioned above, for any $c\geq 0$, the $KG$-module $T^c$ has a submodule isomorphic to $L(c+1)$ with a factor isomorphic to~$(L/L(c+1))^*$. Combining this with Theorem \ref{A13}, we see directly that all $T^c$, $c\geq 0$, have the same composition factors, as predicted by the Brauer-Nesbitt theorem.
\end{note}

\begin{prop} If $0<c<b$ then $soc(T^c)$ contains copies of the non-isomorphic
irreducible modules $M(b)=L(b)$ and $M(0)$. In particular, $soc(T^c)$ is not irreducible.
\end{prop}

\noindent{\it Proof.} The
proof of Theorem \ref{B52} shows that $L(b)$ is inside $T^c$ for all $0\leq c<b$. The map $v\mapsto \ell^c v$
from $I$ into $I(c)$ sends $\ell I$ into $\ell I(c)$, inducing a map from $L$ into $T^c$.
Suppose $c>0$. Let $\l:U\to R^*$ be a group homomorphism such that $P(\l)=G$. Then $\ell^c E_\l$
is in $I(c)$ but not in $I(c+1)$, which shows that the map $L\to T^c$ is not zero.
However, using $c>0$ we easily see that $L(1)$ is in the kernel. Since $L(1)$ is maximal, it follows that
$M(0)=L/L(1)$ embeds into~$T^c$, as claimed.

\begin{note}
 $T^b\cong L^*$, since $T^b=I(b)/\ell I(b)=I(b)/I(b+1)\cong (L/L(b+1))^*=L^*$.
\end{note}

\section{Positive cases of Gow's conjecture}
\label{ultima2}

\begin{thm}
\label{A1}
(a) If $d\leq\ell$ then
(\ref{pato2}) is a composition series of $L$  provided $\lfloor n/e
\rfloor\leq d\ell$.

(b) If $d=\ell+1$ then (\ref{pato2}) is a composition series of $L$
provided $\lfloor n/e \rfloor\leq \ell^2$.

(c) If $d>\ell+1$ then (\ref{pato2}) is a composition series of $L$
provided $\lfloor n/e \rfloor< \ell^2+\ell$.
\end{thm}

\noindent{\it Proof.} This follows from Corollary \ref{icrit}
via Lemma \ref{unodos} of Appendix C.

\begin{note} If $d=1$ Theorem \ref{A1} does not add much to Ackermann's contribution,
as we would just be passing from $\lfloor n/e\rfloor <\ell$ to
$\lfloor n/e\rfloor\leq \ell$. How large can $d$ be? If $\ell|q-1$
and $\ell$ is odd then necessarily $d=1$. However, if $\ell$ is
odd, $2\leq e$ and $e|\ell-1$, or if $\ell=2=e$, then there are
infinitely many primes $q$ such that $q\neq \ell$,
$$e=e(\ell,q)=min\{i\geq 2\,|\, \ell\text{ divides
}\frac{q^i-1}{q-1}\}
$$
and $d>\ell+1$. This follows easily from Dirichlet's Theorem on
primes in arithmetic progression (see Lemma \ref{dri} below for
details). If $q$ is any of these primes then (\ref{filt}) is a
composition series of $L$ as long as $\lfloor n/e\rfloor <
\ell^2+\ell$.
\end{note}

\begin{lem}
\label{dri}
 Let $\ell$ be a prime. If $\ell|q-1$ and $\ell$ is odd then
$d=1$. Suppose that either $\ell=2=e$, or $\ell$ is odd, $2\leq e$
and $e|\ell-1$. Let $s\geq 1$. Then there are infinitely many
primes $q$ such that $q\neq\ell$, $e=e(\ell,q)$ and
$d=\nu_\ell(\frac{q^e-1}{q-1})\geq s$.
\end{lem}

\noindent{\it Proof.} The first assertion follows from the proof
of Lemma \ref{antes} (just replace $es$ by 1).

Suppose still that $\ell$ is odd. Associated to any $m\geq 1$ we
have the multiplicative group~$U(m)=\{[a]\,|\, \gcd(a,m)=1\}$.
Clearly $U(\ell^s)$ decomposes as the direct product of the
kernel, say $A$, of $U(\ell^s)\to U(\ell)$, and a unique subgroup
$B$ isomorphic to~$U(\ell)$. It follows that $U(\ell^s)\to
U(\ell)$ preserves the order of any element whose order divides
$\ell-1$, where all these orders occur since $U(\ell^s)$ is cyclic
of order $(\ell-1)\ell^{s-1}$.

Given $e$ as stated, let $t$ be an integer relatively prime to
$\ell$ having order $e$ modulo $\ell^s$. By Dirichlet's Theorem
there are infinitely many primes congruent to $t$ modulo $\ell^s$.
Let $q$ be one of them. Clearly $q\neq\ell$. The remarks made
above ensure that the order of $q$ modulo $\ell$ is $e$. As $e>1$,
we infer $e=e(\ell,q)$. Moreover, $q^e\equiv t^e\equiv
1\mod\ell^s$, so $d\geq s$.

Suppose next $\ell=2$. By Dirichlet's Theorem there are infinitely
many primes congruent to $-1$ modulo $2^s$, as required.

\begin{lem}
\label{longit} For $i\geq -1$ let $\Lambda_i(n)$ be the total
number of parabolic subgroups of the form
$[z_{-1},z_0,\dots,z_i,0,\dots,0]$ in $\P^*$, as defined in (\ref{demaria}). Then
$\Lambda_{-1}(n)=1$,
$$
\Lambda_i(n)=\underset{0\leq j\leq \lfloor n/e\ell^i\rfloor}\sum
\Lambda_{i-1}(n-e\ell^i j),\quad 0\leq i\leq m,
$$
and $|\P^*|=\Delta_m(n)$.
\end{lem}

\noindent{\it Proof.} This is clear.

\begin{thm}
\label{A9} Suppose the conditions of Theorem \ref{A1}
are satisfied, that is, assume $\lfloor n/e \rfloor\leq d\ell$ if
$d\leq\ell$; $\lfloor n/e \rfloor\leq \ell^2$ if  $d=\ell+1$;
$\lfloor n/e \rfloor< \ell^2+\ell$ if $d>\ell+1$. Then

(a) $c(L)=V=\vert\P^*\vert= \lfloor n/e \rfloor+1$ if $m=0$.

(b) $c(L)=V=\vert\P^*\vert=(x_1+1)(\frac{x_1}{2}\ell+x_0+1)$ if $m=1$ and $\lfloor
n/e\rfloor=(x_1x_0)_\ell$.

(c) $c(L)=V=\vert\P^*\vert=\frac{1}{2}\ell^3+\frac{1}{2}\ell^2+(x_0+1)\ell+2(x_0+1)$ if
$m=2$ and $\lfloor n/e\rfloor=(1,0,x_0)_\ell$.
\end{thm}

\noindent{\it Proof.} We have $c(L)=V$ by Theorem \ref{A1} and $V=\vert\P^*\vert$
by Theorem \ref{muyjota} and Lemma \ref{unodos} of Appendix C. The exact value of
$\vert\P^*\vert$ can be computed by means of Lemma \ref{longit}.

\section{Examples}
\label{lacalle}

\begin{exa}
\label{tatin1}
Suppose
$n=6$, $\ell=2$, $q=5$. In this case, all numbers from 0 to $4=\vartheta(B)$ are $\P$-values.
Thus $V=5$, whereas $|\P^*|=6$. The $3$ bottom factors $M(4)$, $M(3)$, $M(2)$ as well as the top factor $M(0)$ are irreducible. Consider the parabolic subgroups $P=(2,2,2)$ and $Q=(4,1,1)$, where the numbers indicate the sizes of the diagonal blocks. Then $P,Q\in\P^*(1)$, and
James table for $n=6$ adjusted to the prime $\ell=2$ implies that $M(1)=N(P)\oplus N(Q)$ is not irreducible.
\end{exa}

\begin{exa}
\label{tatin2}
Suppose $\ell=2$, $e=2$, $n=10$ and $d=1$ (say $q=5$). Then $V=9$, whereas $|\P^*|=14$. The 14 members of $\P^*$
are distributed into $\P$-values as follows, using an obvious notation for partitions:
$$(82)\in \P^*(0);\,(442), (81^2)\in \P^*(1);\, (4^21^2), (42^3)\in \P^*(2);\,(2^5), (42^21^2)\in \P^*(3);
$$
$$
(2^41^2), (421^4)\in \P^*(4);\, (2^31^4), (41^6)\in\P^*(5);\, (2^21^6)\in\P^*(6);\, (21^8)\in\P^*(7)
;\, (1^{10})\in\P^*(8).
$$
As predicted, the 3 bottom factors $M(8),M(7),M(6)$ as well as the top factor $M(0)$ are irreducible.
Refer now to [3] and use the decomposition matrix from page 257 together with the adjustment matrix from page 258. We see that $c(L)=14$. It follows that all 5 doubtful factors of $L$, namely $M(1)$ through $M(5)$,
fail to be irreducible, and are equal to the direct sum of the two irreducible constituents $N(P),N(Q)$, where
$P$ and $Q$ are as displayed above for each $\P$-value $1\leq c\leq 5$.
\end{exa}

\section{Appendix A}
\label{kmnxz}

The goal of this section is to furnish a proof of Theorem \ref{pringle}.

\subsection{Calculations in the Steinberg lattice}
\label{scsl}

Let $\sigma\in S_n$. The set $I(\sigma)$, of inversions of
$\sigma$, is formed by all pairs $(i,j)$ such that $1\leq i<j\leq
n$ but $\sigma(i)>\sigma(j)$. We associate to $\sigma$ the
subgroup $U_\sigma^+$ formed by all $u\in U$ such that $\sigma
u\sigma^{-1}\in U$, and also the subgroup $U_\sigma^-$ formed by
all $u\in U$ such that $\sigma u\sigma^{-1}\in V$, the lower
unitriangular group. We fix a well-order on $\Phi=\{(i,j)\,|\,
1\leq i<j\leq n\}$. Following this order, we can write any $u\in
U_\sigma^+$ and $v\in U_\sigma^-$ in the form
\begin{equation}
\label{producto} u=\underset{r\notin I(\sigma)}\Pi t_r(a_r) \text{
and } v=\underset{s\in I(\sigma)}\Pi t_s(b_s),
\end{equation}
for unique $a_r,b_s\in F_q$. We have
\begin{equation}
\label{produ2}
U_\sigma^+U_\sigma^-=U=U_\sigma^-U_\sigma^+\quad\text{ and }\quad
U_\sigma^+\cap U_\sigma^-=1.
\end{equation}
For the special permutation
\begin{equation}
\label{design} \sigma_0=(1,n)(2,n-1)(3,n-2)\cdots=\sigma_0^{-1}
\end{equation}
we have $I(\sigma_0)=\Phi$, so that
\begin{equation}
\label{produ3} U_{\sigma_0}^-=U\quad\text{ and }\quad
U_{\sigma_0}^+=1.
\end{equation}
Moreover,
\begin{equation}
\label{produ4} I(\sigma_0\sigma)=\Phi\setminus
I(\sigma)\quad\text{ and }\quad U_{\sigma_0\sigma}^+=U_\sigma^-.
\end{equation}

The subset $\{g\widehat{B}\,|\, g\in G\}$ of $RG$ is linearly
independent, so it is an $R$-basis for its span, say $Y$. Note
that $I$ is contained in $Y$. If $x\in I$ it is then clear what we
mean by ``the coefficient of $g\widehat{B}$ in $x$", a phrase that
will be used at critical points below. Of course, we  may have
$g\widehat{B}=h\widehat{B}$ for $g,h\in G$, which happens if and
only if $gB=hB$. We can avoid repetitions by means of the Bruhat
decomposition. Thus, a basis for $Y$ is formed by all $u\sigma
\widehat{B}$, where $\sigma\in S_n$ and $u\in U_{\sigma^{-1}}^-$.

The following two results are valid in the more general context
used in [4].

\begin{lem}
\label{ayuda}
 Let $\l:U\to R^*$ be a group homomorphism with $E_\l$ as in (\ref{defegrande}). Then
\begin{equation}
\label{suma} E_\l=\underset{\sigma\in S_n}\sum\; \underset{u\in
U_{\sigma^{-1}}^-}\sum sg(\sigma) C_\sigma(\l) \l(u) u\sigma
\widehat{B},
\end{equation}
where
\begin{equation}
\label{arbol} C_\sigma(\l)= \underset{v\in U_{\sigma^{-1}}^+}\sum
\l(v)=\begin{cases} |U_{\sigma^{-1}}^+| & \text{ if }\l\text{
is trivial on }U_{\sigma^{-1}}^+,\\
0 & \text{ otherwise }.
\end{cases}
\end{equation}
\end{lem}

\noindent{\it Proof.} According to the definitions
(\ref{defechico}) of $\u$ and (\ref{defegrande}) of $E_\l$ we have
$$
E_\l= \underset{u\in U}\sum \l(u)u\underset{\sigma\in S_n}\sum
sg(\sigma)\sigma \widehat{B}=\underset{\sigma\in S_n}\sum
\underset{u\in U}\sum sg(\sigma)\l(u)u\sigma \widehat{B}.
$$
We now use the decomposition (\ref{produ2}) of $U$, the fact that
$\sigma^{-1}v\sigma\widehat{B}=\widehat{B}$ for all $v\in
U_{\sigma^{-1}}^+$, and that $\l$ is a group homomorphism to
obtain (\ref{suma}). The displayed value of $C_\sigma(\l)$ is
clear.

\begin{lem}
\label{ayu2} Let $\sigma\in S_n$. Let $\l,\mu:U\to R^*$ be group
homomorphisms. Suppose that every $X_r$, $r\in\Pi$, acts on the
element $\widehat{U_{\sigma^{-1}}^-}\cdot \sigma \cdot E_\l$ of
$I$ via $\mu^{-1}$. Then
$$
\widehat{U_{\sigma^{-1}}^-}\cdot \sigma \cdot
E_\l=sg(\sigma)E_\mu.
$$
\end{lem}

\noindent{\it Proof.} Since the $X_r$, $r\in\Pi$, generate  $U$,
it follows that $U$ acts on $\widehat{U_{\sigma^{-1}}^-}\sigma
E_\l$ via $\mu^{-1}$. But $U$ acts on $I$ via the regular
representation. We deduce that $\widehat{U_{\sigma^{-1}}^-}\sigma
E_\l$ must be a scalar multiple of $E_\mu$, that is
\begin{equation}
\label{cai} \widehat{U_{\sigma^{-1}}^-}\sigma E_\l=a E_\mu,
\end{equation}
where $a\in R$ is to be found.  To determine $a$ we write both
sides of (\ref{cai}) relative to the basis $\{g\widehat{B}\,|\,
g\in G\}$ of $Y$ previously mentioned, and compare coefficients.
In view of (\ref{cai}), it suffices to compare coefficients in a
{\em single} basis vector $g\widehat{B}$, provided the coefficient
of $g\widehat{B}$ in $E_\mu$ is not zero. A good choice turns out
to be $\sigma\sigma_0\widehat{B}$, where $\sigma_0$ is defined in
(\ref{design}).

By (\ref{suma}) and (\ref{produ3}), the coefficient of
$\sigma_0\widehat{B}$ in $E_\l$ is equal to $sg(\sigma_0)$.
Multiplication  by $\sigma$ simply shifts all basis vectors,
so the coefficient of $\sigma\sigma_0\widehat{B}$ in $\sigma E_\l$
is also $sg(\sigma_0)$.

Now by (\ref{produ4})
$$
U_{(\sigma\sigma_0)^{-1}}^+=U_{\sigma_0^{-1}\sigma^{-1}}^+=
U_{\sigma_0\sigma^{-1}}^+=U_{\sigma^{-1}}^-.
$$
Thus if $u\in U_{\sigma^{-1}}^-=U_{(\sigma\sigma_0)^{-1}}^+$ then
$$
u \sigma\sigma_0 \widehat{B}= \sigma\sigma_0[(\sigma\sigma_0)^{-1}
u \sigma\sigma_0 ]\widehat{B}=\sigma\sigma_0 \widehat{B},
$$
so multiplying $\sigma E_\l$ by $u$ fixes the basis vector
$\sigma\sigma_0 \widehat{B}$. This happens for the
$|U_{\sigma^{-1}}^-|$ vectors $u$ in $U_{\sigma^{-1}}^-$, which,
so far, will produce the coefficient
$sg(\sigma_0)|U_{\sigma^{-1}}^-|$ for $\sigma\sigma_0 \widehat{B}$
in $\widehat{U_{\sigma^{-1}}^-}\sigma E_\l$.

We must now make sure that the basis vector $\sigma\sigma_0
\widehat{B}$ cannot be produced in any other way in
$\widehat{U_{\sigma^{-1}}^-}\sigma E_\l$. Well, by (\ref{suma}), a
typical summand of $E_\l$ has the form $v\tau \widehat{B}$, where
$\tau\in S_n$ and $v\in U_{\tau^{-1}}^-$.  Thus, a typical summand
of $\widehat{U_{\sigma^{-1}}^-}\sigma E_\l$ will have the form
$u\sigma v\tau \widehat{B}$, where $u\in U_{\sigma^{-1}}^-$.  When
will this summand equal $\sigma\sigma_0 \widehat{B}$? Well,
suppose that $u\sigma v\tau \widehat{B}=\sigma\sigma_0
\widehat{B}$ for some $u$, $v$ and $\tau$ as stated. The right
hand side was shown above to equal $u\sigma\sigma_0 \widehat{B}$,
which gives $u\sigma v\tau \widehat{B}=u\sigma\sigma_0
\widehat{B}$, and a fortiori the equation $u\sigma v\tau
B=u\sigma\sigma_0 B$ in $G$. This, in turn, yields $v\tau
B=\sigma_0 B$. The uniqueness part of the Bruhat decomposition
gives $\tau=\sigma_0$ first, and then $v=1$, since
$U_{\sigma_0}^-=U$. Thus, the basis vector $\sigma\sigma_0
\widehat{B}$ appears in $\widehat{U_{\sigma^{-1}}^-}\sigma E_\l$
only as described above. Hence the coefficient of $\sigma\sigma_0
\widehat{B}$ in $\widehat{U_{\sigma^{-1}}^-}\sigma E_\l$ is
exactly $sg(\sigma_0)|U_{\sigma^{-1}}^-|$. In particular,
$\widehat{U_{\sigma^{-1}}^-}\sigma E_\l\neq 0$.

Observe next that $\mu$ is trivial on $U_{\sigma^{-1}}^-$. Indeed,
let $u\in U_{\sigma^{-1}}^-$. Clearly, $u
\widehat{U_{\sigma^{-1}}^-}\sigma
E_\l=\widehat{U_{\sigma^{-1}}^-}\sigma E_\l$, while by hypothesis
$u \widehat{U_{\sigma^{-1}}^-}\sigma
E_\l=\mu(u)^{-1}\widehat{U_{\sigma^{-1}}^-}\sigma E_\l$. Since $I$
is a torsion free $R$-module and, by above,
$\widehat{U_{\sigma^{-1}}^-}\sigma E_\l\neq 0$, we infer
$\mu(u)=1$.

Finally, due to (\ref{suma}), the coefficient of
$\sigma\sigma_0\widehat{B}$ in $E_\mu$ is equal to
$sg(\sigma\sigma_0)C_{\sigma\sigma_0}(\mu)$. By above $\mu$ is
trivial on $U_{\sigma^{-1}}^-=U_{(\sigma\sigma_0)^{-1}}^+$.
Therefore (\ref{arbol}) gives
$C_{\sigma\sigma_0}(\mu)=|U_{\sigma^{-1}}^-|$. Hence the
coefficient of $\sigma\sigma_0\widehat{B}$ in $E_\mu$ is equal to
$sg(\sigma)sg(\sigma_0)|U_{\sigma^{-1}}^-|$.

Comparing coefficients yields $a=sg(\sigma)$, as claimed.

\subsection{Properties of parabolic subgroups reflected on $I$}
\label{huija}

Let $P=(a_1,...,a_k)$ be a parabolic subgroup. Replacing any $a_i>1$ by a subsequence $(a,b)$ such that $a+b=a_i$
produces a parabolic subgroup contained in $P$, and any parabolic
subgroup contained in $P$ can be obtained by repeated application
of this procedure.

Let $J$ be the subset of $\Pi$ corresponding to $P$. It is clear
what we mean by the connected components of $J$. We next describe
how these can be read off from $(a_1,...,a_k)$. If $a_1=1$ then
$(1,2)$ is not in $J$, while if $a_1>1$ then all of
$(1,2),...,(a_1-1,a_1)$ are in $J$ but $(a_1,a_1+1)$ is not in
$J$. The same procedure is applied to $a_2,...,a_k$, starting at
the first element of $\Pi$ whose inclusion in $J$ was not decided
in the previous steps. For instance, $P=(2,1,2)$ produces
$J=\{(1,2),(4,5)\}$. Each $a_i>1$ gives rise to a connected
component of $J$ of length $a_i-1$, and every connected component
of $J$ arises in this way. Let $Q$ be the parabolic subgroup
obtained from $J$ by a single switching $a_i\leftrightarrow
a_{i+1}$. Let $J'$ be the subset of $\Pi$ associated to $Q$. How
is $J'$ obtained from $J$? This is obvious, but later applications
of Lemma \ref{ayu2} will require an explicit answer. Four cases
arise:

$\bullet$ Suppose $a_i=a_{i+1}=1$. Then $J'=J$.

$\bullet$ Suppose $a_i>1$ and $a_{i+1}>1$. Let
$$
A=\{(j,j+1),...,(j+m-1,j+m)\},\quad m\geq 1
$$
and
$$
B=\{(j+m+1,j+m+2),...,(j+m+s,j+m+s+1)\},\quad s\geq 1
$$
be the connected components of $J$ corresponding to $a_i=m+1$ and
$a_{i+1}=s+1$. Then the connected components of $J'$ are precisely
those of $J$, except for $A$, which must be replaced by
$$
A'=\{(j,j+1),...,(j+s-1,j+s)\},
$$
and for $B$, which must be replaced by
$$
B'=\{(j+s+1,j+s+2),...,(j+s+m,j+s+m+1)\}.
$$
Of course, $J'=J$ if $a_i=a_{i+1}$. Note that $(j+m,j+m+1)\notin
J$, while $(j+s,j+s+1)\notin J'$.

$\bullet$ Suppose $a_i>1$ and $a_{i+1}=1$. Then $a_i=m+1$, where
$m\geq 1$. Denote by $A=\{(j,j+1),...,(j+m-1,j+m)\}$ the connected
component of $J$ associated to $a_i$. In this case $J'$ has the
same connected components as $J$, except for $A$, which must be
replaced by $A'=\{(j+1,j+2),...,(j+m,j+m+1)\}$.

$\bullet$ Suppose $a_i=1$ and $a_{i+1}>1$. Then $a_{i+1}=s+1$,
where $s\geq 1$. Denote by $A=\{(j+1,j+2),...,(j+s,j+s+1)\}$
the connected component of $J$ associated to $a_{i+1}$. In this
case $J'$ has the same connected components as $J$, except for
$A$, which must be replaced by $A'=\{(j,j+1),...,(j+s-1,j+s)\}$.

\begin{thm}
\label{barb1} If $P,Q\in\P$ are equivalent then $I'(P)=I'(Q)$.
\end{thm}

\noindent{\it Proof.} Let $P=(a_1,...,a_k)$ and let $J$ be the
subset of $\Pi$ associated to $P$. It suffices to prove the
theorem when $Q$ is obtained from $P$ by a single switching
$a_i\leftrightarrow a_{i+1}$. Let $J'$ be the subset of $\Pi$
associated to $Q$.

Our main tool will be Lemma \ref{ayu2}. Once the right choice of
$\sigma\in S_n$ is made, it is then a matter of routine to verify
that the hypotheses of Lemma \ref{ayu2} are met.

We refer to the notation introduced earlier in this section for
this scenario. Of the four given cases, we only need to consider
the last three. Let us begin with the first of these, namely when
$a_i>1$ and $a_{i+1}>1$.

Let $\sigma\in S_n$ fix every point outside of the interval
$[j,...,j+m+s+1]$ and be defined as follows on this interval:
$$
\begin{array}{cccccc}
  j & \cdots & j+m & j+m+1 & \cdots & j+m+s+1 \\
\downarrow & \cdots & \downarrow & \downarrow & \cdots & \downarrow\\
  j+s+1 & \cdots & j+m+s+1 & j & \cdots & j+s \\
\end{array}
$$
Notice that
$$
\sigma A\sigma^{-1}=B'\text{ and }\sigma B\sigma^{-1}=A'.
$$
Thus $\sigma J\sigma^{-1}=J'$ and conjugation by $\sigma$ sends
the connected components of $J$ into those of $J'$.

Clearly conjugation by the non-trivial permutation $\sigma$ cannot
preserve $\Pi$. In this case, the following subsets of $\Pi$ are
sent outside of $\Pi$: the ``middle" set $C=\{(j+m,j+m+1)\}$ and
the ``boundary" set $D=\{(j-1,j),(j+m+s+1,j+m+s+2)\}\cap \Pi$.
Also notice that conjugation by $\sigma$ does not send $P$ into
$Q$ either. Indeed, if $s\neq m$ then $P\neq Q$, and distinct
standard parabolic subgroups cannot be conjugate, while if $s=m$
then $P=Q$, but still $\sigma\notin P$, and $P$ is
self-normalizing.

Let $\l:U\to R^*$ be a group homomorphism such that $P(\l)=P$. We
next define a group homomorphism $\mu:U\to R^*$ such that
$P(\mu)=Q$. It suffices to define a group homomorphism on every
$X_r$, $r\in \Pi$, as these will have a unique extension to $U$
(we use here that there are exactly $|U/U'|$ homomorphisms $U\to
R^*$, given that $U/U'$ is an elementary abelian $p$-group and $R$
has a non-trivial $p$-root of unity). We simply let
\begin{equation}
\label{mu1} \mu(t_r(a))=\l(t_{\sigma^{-1}r\sigma}(a)),\quad r\in
J'
\end{equation}
and
\begin{equation}
\label{mu2}
\mu(t_r(a))=1,\quad r\in\Pi\setminus J'.
\end{equation}
By construction, $P(\mu)=Q$.

By virtue of Lemma \ref{ayu2}, all we have to do now is verify
that each fundamental root subgroup acts on
$\widehat{U_{\sigma^{-1}}^-}\cdot \sigma \cdot E_\l$ via
$\mu^{-1}$. Indeed, this will show that $I'(Q)\subseteq I'(P)$,
and switching back $a_i$ and $a_{i+1}$ will yield the reverse
inclusion.

Note first of all that
\begin{equation}
\label{invers} I(\sigma^{-1})=\{(a,b)\,|\quad j\leq a\leq
j+s,\quad j+s+1\leq b\leq j+m+s+1\}.
\end{equation}

We next verify that each $X_r$, $r\in\Pi$, acts on
$\widehat{U_{\sigma^{-1}}^-}\cdot \sigma \cdot E_\l$ via
$\mu^{-1}$. Now $\Pi$ decomposes as $\Pi=A'\cup B'\cup C'\cup
D\cup E$, where $C'=\{(j+s,j+s+1)\}$, $A'$, $B'$ and $D$ have been defined above,
and $E$ is the complement of $A'\cup B'\cup C'\cup D$ in $\Pi$.
Our argument is divided according to this decomposition.

If $r$ is in $E$ then $X_r$ normalizes $U_{\sigma^{-1}}^-$ and
commutes elementwise with $\sigma$, so it acts on
$U_{\sigma^{-1}}^-\sigma E_\l$ via $\l^{-1}$, and hence via
$\mu^{-1}$, as they agree on $X_r$.

If $r=(j+s,j+s+1)$ then $X_r$ is included in $U_{\sigma^{-1}}^-$,
so it acts trivially on $\widehat{U_{\sigma^{-1}}^-}\sigma E_\l$,
and hence via $\mu^{-1}$, since, as remarked earlier, $(j+s,j+s+1)\notin J'$.

Consider next the case when $r\in A'\cup B'$.  We will make use of
the well-known formula:
\begin{equation}
\label{conjugar} \sigma
t_{ij}(a)\sigma^{-1}=t_{\sigma(i)\sigma(j)}(a),\quad \sigma\in
S_n.
\end{equation}
We will also use the commutator $[xy]=xyx^{-1}y^{-1}$. Clearly if
$i<j$, $k<l$ and $i\neq l$ then
\begin{equation}
\label{commutator} [t_{ij}(a)t_{kl}(b)]=\begin{cases} t_{il}(ab) &
\text{ if } j=k,\\ 1 & \text{ otherwise }. \end{cases}
\end{equation}
From (\ref{commutator}) and (\ref{invers}) we see that $X_r$
normalizes $U_{\sigma^{-1}}^-$. Thus by (\ref{conjugar})
$$
t_r(a)\widehat{U_{\sigma^{-1}}^-}\sigma
E_\l=\widehat{U_{\sigma^{-1}}^-}t_r(a)\sigma
E_\l=\widehat{U_{\sigma^{-1}}^-}\sigma\sigma^{-1}t_r(a)\sigma
E_\l=\widehat{U_{\sigma^{-1}}^-}\sigma
t_{\sigma^{-1}r\sigma}(a)E_\l,
$$
where the last term equals
$$
\l(t_{\sigma^{-1}r\sigma})^{-1}\widehat{U_{\sigma^{-1}}^-}\sigma
E_\l=\mu(t_r(a))^{-1}\widehat{U_{\sigma^{-1}}^-}\sigma E_\l.
$$

Suppose finally that $r$ belongs to $D$. Let us treat the case
$r=(j-1,j)$ first. It is no longer true that $X_r$ normalizes
$U_{\sigma^{-1}}^-$, so we have to be a bit careful. Let
$t_r(\alpha)\in X_r$ and let $u=U_{\sigma^{-1}}^-$. Selecting a
suitable ordering, we may use (\ref{producto}) to write
$u=u_1u_2$, where $u_1$ is a product of factors of the form
$t_{ab}(\beta)$, where $(a,b)\in I(\sigma^{-1})$ and $a\neq j$,
and $u_2$ is a product of factors of the form $t_{jb}(\beta)$,
where $(j,b)\in I(\sigma^{-1})$. By (\ref{commutator}) we have
$$
t_r(\alpha)u_1=u_1t_r(\alpha).
$$
By (\ref{commutator}) any $t_{jb}(\beta)$ will commute with any
commutator
$$
[t_r(\alpha)t_{jc}(\gamma)]=t_{j-1,c}(\delta),
$$
where $j+s+1\leq b,c\leq j+m+s+1$. Repeatedly using this comment and the given expression for $u_2$, we see
that $t_r(\alpha)u_2=u_2t_r(\alpha)z$, where $z$ is a product of
factors of form $t_{j-1,c}(\delta)$, where $j+s+1\leq c\leq
j+m+s+1$. Therefore $t_r(\alpha)u=ut_r(\alpha)z$. Now
$w=\sigma^{-1}z\sigma$ is a product of factors of the form
$t_{j-1,d}(\delta)$, where $j\leq d\leq j+m$. Now if $d>j$ then
$t_{j-1,d}(\delta)\in U'$, while $t_{j-1,j}(\delta)$ acts
trivially on $E_\l$, since $(j-1,j)\notin J$. Thus $w$ acts
trivially on $E_\l$. Also
$\sigma^{-1}t_{j-1,j}(\alpha)\sigma=t_{j-1,j+m+1}(\alpha)\in U'$
acts trivially on $E_\l$. All in all, we get that $t_r(\alpha)$
acts trivially on $u\sigma E_\l$. As this happens for all $u\in
U_{\sigma^{-1}}^-$, we finally obtain that $t_r(\alpha)$ acts
trivially on $\widehat{U_{\sigma^{-1}}^-}\sigma E_\l$. The
reasoning when $r=(j+m+s+1,j+m+s+2)$ is entirely analogous.

This completes the proof of the case $a_i>1$ and $a_{i+1}>1$. The
case $a_i>1$ and $a_{i+1}=1$ can be handled as a degenerate (and
simplified) case of the above, corresponding to $s=0$.
Accordingly, we merely need to modify the permutation $\sigma$ to
$$
\begin{array}{ccccc}
  j & j+1 & \cdots & j+m & j+m+1 \\
\downarrow & \downarrow & \cdots & \downarrow & \downarrow\\
  j+1 & j+2 & \cdots & j+m+1 & j \\
\end{array}
$$
Similarly, the case $a_i=1$ and $a_{i+1}>1$ can also be handled as
a degenerate case of the one above, corresponding to $m=0$. Here
we modify $\sigma$ to the permutation
$$
\begin{array}{ccccc}
  j & j+1 & \cdots & j+s & j+s+1 \\
\downarrow & \downarrow  & \cdots & \downarrow & \downarrow\\
  j+s+1 & j & \cdots & j+s-1 & j+s \\
\end{array}
$$
In the notation corresponding to these cases, conjugation by
$\sigma$ will send $A$ to $A'$ and fix all other connected
components of $J$. Given a group homomorphism $\l:U\to R^*$ such
that $P=P(\l)$, we define $\mu$ using the formulae (\ref{mu1}) and
(\ref{mu2}). Again, $P(\mu)=Q$, and one can check that the
argument given in the general case will go through in the two
degenerate cases above, mutatis mutandi.

\begin{thm}
\label{barb2}
 Let $Q\subseteq P$ be parabolic subgroups of $G$. Then
$I'(Q)\subseteq I'(P)$.
\end{thm}

\noindent{\it Proof.} Let $J$ and $J'$ be the subsets of $\Pi$
associated to $P$ and $Q$, respectively. We may assume that
$J\neq\emptyset$ and $J'\neq J$. By repeatedly removing one point
from $J$ at a time, we may assume that $J'$ is obtained by
removing a single point, say $r$, from $J$. Thus
$J'=J\setminus\{r\}$. Let $A$ be the connected component of $J$ to
which $r$ belongs. Two cases arise: $r$ is an endpoint or $r$ is a
middle point of $A$.

Now an endpoint can be a left or a right endpoint. A middle point
can be skewed to the left, i.e. there are at least as many points
in $A$ to the right of it as to the left of it, or skewed to the
right. By means to Theorem \ref{barb1} we may reduce ourselves to
consider only left endpoints and middle points skewed to the left.

This is so because the bijection $(1,2)\leftrightarrow (n-1,n),
(2,3)\leftrightarrow (n-2,n-1),\dots$ of $\Pi$ into itself induces
a bijection from $\P$ into itself, which sends a parabolic
subgroup into one equivalent to it, and interchanges left and
right in both cases above.

By rearranging the blocks of $P$ and using Theorem \ref{barb1}, we
may also assume that the left endpoint of $A$ is $(1,2)$. Thus
$A=\{(1,2),...,(k-1,k)\}$, where $k>1$.

Assume first that $r$ is the left endpoint of $A$, so that
$r=(1,2)$. Then $J'$ has the same connected components as $J$,
except for $A$, which must now be replaced by
$A'=\{(2,3),...,(k-1,k)\}$. Note that $A=\emptyset$ if $k=2$.

Consider the cycle $\sigma=(1,2,...,k)\in S_n$. Given a group
homomorphism $\l:U\to R^*$ such that $P(\l)=P$, we define $\mu$
using (\ref{mu1}) and (\ref{mu2}). Then $P(\mu)=Q$. We now apply
Lemma \ref{ayu2}, verifying its hypotheses as in the proof Theorem
\ref{barb1}.

Suppose next $r=(i,i+1)$ is a middle point of $A$ skewed to the
left. Thus
$$
A=\{(1,2),...,(i-1,i),(i,i+1),(i+1,i+2),...,(2i-1,2i),...,(k-1,k)\},
$$
where $1<i$ and $2i\leq k$. The connected components of $J'$ are
those of $J$, except that $A$ must be replaced by the two
components
$$
A'=\{(1,2),...,(i-1,i)\}\text{ and
}B=\{(i+1,i+2),...,(2i-1,2i),...,(k-1,k)\}.
$$
Consider the permutation $\sigma\in S_n$ whose inverse
$\sigma^{-1}$ fixes every number larger than $k$ and has the
following effect on the interval $[1,...,k]$:
$$
\begin{array}{ccccccccccc}
  i+1 & i+2 & \cdots & 2i-1 & 2i & \cdots & k & 1 & 2 & \cdots & i \\
  \downarrow & \downarrow & \cdots & \downarrow & \downarrow &
  \cdots & \downarrow & \downarrow & \downarrow &
  \cdots & \downarrow \\
  1 & 2 & \cdots & i-1 & i & \cdots & k-i & k-i+1 & k-i+2 & \cdots & k \\
\end{array}
$$
This definition of $\sigma^{-1}$ yields
$$
I(\sigma^{-1})=\{(a,b)\,|\, 1\leq a\leq i,\quad i+1\leq b\leq k\}.
$$
As usual, a valid application of Lemma \ref{ayu2} yields the
desired result.

\begin{note} Various special cases suggest that $[P:Q]I'(P)\subseteq I'(Q)$
if $Q\subseteq P$ are in $\P$.
\end{note}

\section{Appendix B}
\label{snumco}

Here we develop auxiliary tools to compute $\nu_\ell([P:B])$. Recall that $d$ is defined in (\ref{dedd}) and
$$
s_0=d, s_1=\ell d+1, s_2=\ell^2d+\ell+1, s_3=\ell^3d+\ell^2+\ell+1,\dots.
$$
For typographical reasons it will be necessary to use the notation
$$
w(a,b)=\frac{q^a-1}{q^b-1},\quad w(a)=\frac{q^a-1}{q-1},\quad
a,b\geq 1,
$$
$$
g(a)=\nu_\ell(w(a)),\quad h(a)=\nu_\ell(w(1)w(2)\cdots w(a)),\quad
a\geq 1.
$$
The following two results are borrowed from [2].

\begin{lem}
\label{antes} Let $s$ be a positive integer. Then
$$
\nu_\ell\big[\frac{q^{es\ell}-1}{q^{es}-1}\big]=1.
$$
\end{lem}
\noindent{\it Proof.} Suppose first that $\ell=2$. Then $e=2$, $q$
is odd and
$$
\frac{q^{es2}-1}{q^{es}-1}=q^{es}+1=(q^s)^2+1\equiv 2\mod 4.
$$
Suppose next $\ell>2$. We have $q^{es}-1=a\ell^b$, with $a$
coprime to $\ell$ and $b\geq 1$. Then
$$
\frac{q^{es\ell}-1}{q^{es}-1}=\frac{(a\ell^b+1)^\ell-1}{a\ell^b}=\underset{1\leq
i\leq \ell}\sum{\ell\choose i} (a\ell^b)^{i-1}\equiv
\ell\mod\ell^2.
$$

\begin{lem}
\label{james} Let $t$ be a positive integer. Then
$$
\nu_\ell\big[\frac{q^{et}-1}{q^e-1}\big]=\nu_\ell(t).
$$
\end{lem}

\noindent{\it Proof.} We have $t=c\ell^u$, with $c$ coprime to
$\ell$. Then
$$
\frac{q^{et}-1}{q^e-1}=\frac{q^{ec}-1}{q^e-1}\times\underset{1\leq
i\leq u}\Pi w(ec\ell^i,ec\ell^{i-1}).
$$
But
$$
\frac{q^{ec}-1}{q^e-1}\equiv 1+q^e+\cdots+q^{e(c-1)}\equiv
c\not\equiv 0\mod \ell,
$$
while $\ell$ divides each factor $w(ec\ell^i,ec\ell^{i-1})$
 exactly once by Lemma
\ref{antes}, so the result follows.

\begin{lem}
\label{v8} $h(e\ell^i)=s_i$ for all $i\geq 0$.
\end{lem}

\noindent{\it Proof.} First note that by Lemma \ref{james}
\begin{equation}
\label{esete} \nu_\ell(s)=\nu_\ell(t)\Rightarrow g(es)=g(et),\quad
s,t\geq 1.
\end{equation}
Next observe that $\ell\mid w(a)$ if and only if $e|a$. It follows
from this observation that if $a=be+c$, where $0\leq b$ and $0\leq
c<e$, then
\begin{equation}
\label{xxx} h(a)=\underset{1\leq i\leq b}\sum g(ie)=h(be).
\end{equation}
We deduce from (\ref{xxx}) that $h(e)=g(e)=d$, so our formula
works if $i=0$. Suppose $h(e\ell^i)=s_i$ for some $i\geq 0$. Then
by (\ref{xxx})
$$
h(e\ell^{i+1})=h(e\ell^i)+\underset{1\leq k\leq \ell^i}\sum
g(e(k+\ell^i))+\cdots+\underset{1\leq k\leq \ell^i}\sum
g(e(k+(\ell-1)\ell^i)).
$$
If $1\leq k\leq \ell^i$ and $0\leq j<\ell-1$, or if $1\leq k<
\ell^i$ and $j=\ell-1$, then $\nu_\ell(k+j\ell^i)=\nu_\ell(k)$. On
the other hand if $k=\ell^i$ and $j=\ell-1$ then
$\nu_\ell(k+j\ell^i)=\nu_\ell(\ell^{i+1})=\nu_\ell(\ell^i)+1$. We
infer from (\ref{esete}) that
$$
h(e\ell^{i+1})=\underbrace{h(e\ell^i)+h(e\ell^i)+\cdots+h(e\ell^i)}_{\ell}+1=\ell
s_i+1=s_{i+1}.
$$

\begin{lem}
\label{u1}

Let $a=ex$, where $x=b\ell^i+y$, $0\leq i$, $0\leq b<\ell$ and
$0\leq y<\ell^i$. Then
$$
h(a)=bh(e\ell^i)+h(ey).
$$
\end{lem}

\noindent{\it Proof.} We have $a=eb\ell^i+ey$, where by
(\ref{xxx})
$$
h(a)=h(e\ell^i)+\underset{1\leq k\leq \ell^i}\sum
g(e(k+\ell^i))+\cdots+\underset{1\leq k\leq \ell^i}\sum
g(e(k+(b-1)\ell^i))+\underset{1\leq k\leq y}\sum g(e(k+b\ell^i)).
$$
If $1\leq k\leq \ell^i$ and $0\leq j<b$, or if $1\leq k<\ell^i$
and $j=b$, then $\nu_\ell(k+j\ell^i)=\nu_\ell(k)$. The rest
follows much as above.

\begin{lem}
\label{r4} If $1\leq a\leq n$ and $\Delta(a)=(y_{-1},y_0,\dots,y_m)$
is as given in Section \ref{cmp} then
$$
h(a)=y_mh(e\ell^m)+\cdots+y_1h(e\ell)+y_0h(e)=y_m
s_m+\cdots+y_1s_1+y_0s_0.
$$
\end{lem}

\noindent{\it Proof.} This follows by using Lemmas \ref{v8} and
\ref{u1}, as well as (\ref{xxx}).

\section{Appendix C}
\label{cxcb}

Here we determine when the map $P\mapsto \nu_\ell([G:P])$ is injective on $P^*$, i.e.
when $\vert \P^*(c)\vert=1$ for all $\P$-values $c$. We adopt all of the notation
introduced in Section \ref{noahuant}. Clearly, the injectivity of $\t$ on $\P^*$ is equivalent
to the injectivity of $\p$ on $\P^*$.

We define $\HP$ to be the set of all $[z_0,z_1,0,\dots,0]\in\P^*$. Note that
$\HP=\P^*$ if $\lfloor n/e\rfloor<\ell^2$.

\begin{lem}
\label{zx1} $\p$ is injective on $\HP$ if and only if $\lfloor
n/e\rfloor\le d\e$.
\end{lem}

\noindent{\it Proof.} Suppose $\lfloor n/e\rfloor\le d\e$ and
$\p([z_0,z_1,0,\dots,0])=\p([z_0',z_1',0,\dots,0])$. By Theorem
\ref{muyjota}
$$
z_0d+z_1(d\ell+1)= z_0'd+z_1'(d\ell+1).
$$
Since $\gcd(d,d\e+1)=1$ there must be an integer $k$ such that
\begin{equation}
\label{cuak} (z_0',z_1')=(z_0+k(d\ell+1),z_1-kd).
\end{equation}
Now $\lfloor n/e\rfloor\le d\e$ forces $0\leq z_0,z_0'\le d\e$, so
(\ref{cuak}) implies $z_0'=z_0$, and a fortiori $z_1'=z_1$.

Suppose next $\lfloor n/e\rfloor\geq d\e+1$. Then
$P=[d\e+1,0,\dots,0], Q=[0,d,0,\dots,0]\in \HP$ and
$$
\p(P)=(d\e+1)d=\p(Q)
$$
by Theorem \ref{muyjota}, so $\p$ is not injective on $\HP$.

\begin{lem}
\label{zx4} If $\lfloor n/e\rfloor\geq \ell^2+\ell$ then $\p$ is
not injective on $\P^*$.
\end{lem}

\noindent{\it Proof.} Let $P=[\ell,0,1,0,\dots,0]$ and
$Q=[0,\ell+1,0,\dots,0]$. Then $P,Q\in\P^*$ and
$$
\p(P)=d\ell+\ell(d\ell+1)+1=d\ell(1+\ell)+1+\ell=(1+\ell)(d\ell+1)=\p(Q).
$$

\begin{lem}
\label{zx3} Suppose that $\e^2\leq \lfloor n/e\rfloor
<\ell^2+\ell$ and $\lfloor n/e\rfloor\leq \ell d$. Then there are
no parabolic subgroups $P=[z_0,z_1,0,\dots,0]\in\HP$ and
$Q=[a,0,1,0,\dots,0]\in\P^*$ such that $\p(P)=\p(Q)$, except only
if $d=\ell+1$ and $\ell^2+1\leq  \lfloor n/e\rfloor$, when such
$P$ and $Q$ do exist.
\end{lem}

\noindent{\it Proof.} Suppose $\p(P)=\p(Q)$ for $P,Q$ as stated.
Then by Theorem \ref{muyjota}
\begin{equation}
\label{hy1} z_0d+z_1(d\e+1)=ad+\e(d\e+1)+1.
\end{equation}
Hence there is an integer $k$ such that
\begin{equation}
\label{kint} z_0=a-\ell+k(d\ell+1),\quad z_1=1+\ell-kd.
\end{equation}
If $k\leq 0$ then $1+\ell-dk \geq 1+\ell$, against the fact that
$\lfloor n/e\rfloor<\ell(\ell+1)$. Therefore $k>0$.

Observe now that our hypotheses imply $\ell\leq d$. If $k\geq 3$
then $1+\ell-kd<0$, which is impossible. If $k=2$ then
$1+\ell-2d\geq 0$ implies $d=1=\ell$, which is absurd. The only
possibility is $k=1$ with $d=\ell$ or $d=\ell+1$.

If $d=\ell$ our hypotheses yield $\lfloor n/e\rfloor=\ell^2$. Then
from $Q=[a,0,1,0,\dots,0]\in\P^*$ we infer $a=0$. Replacing the
values $k=1$, $a=0$ and $d=\ell$ in (\ref{kint}) gives
$z_0=\ell^2-\ell+1$ and $z_1=1$. Then
$z_0+z_1\ell=\e^2-\e+1+\e=\e^2+1$, contradicting the fact that
$\lfloor n/e\rfloor=\ell^2$.

All in all, we must have $k=1$ and $d=\ell+1$. Going back to
(\ref{kint}) we obtain $z_1=0$ and $z_0=\ell^2+1+a$. In particular
$\lfloor n/e\rfloor\geq \ell^2+1$. This shows that no such $P,Q$
exist, except when $d=\ell+1$ and $\ell^2+1\leq \lfloor
n/e\rfloor< \ell^2+\ell$. In this last case, setting $k=1$, $a=0$
yields $P=[\ell^2+1,0,\dots,0]\in\HP$ and
$Q=[0,0,1,0,\dots,0]\in\P$, with $\p(P)=(\ell+1)(\ell^2+1)=\p(Q)$.
The simplest example occurs when $\e=2$, $q=7$ and $n=10$.

\begin{lem}
\label{unodos}
 (a) Suppose $d\leq \ell$. Then $\t$
is injective on $\P^*$ if and only if $\lfloor n/e\rfloor\leq
d\ell$.

(b) Suppose $d=\ell+1$. Then $\t$ is injective on $\P^*$ if and
only if $\lfloor n/e\rfloor\leq \ell^2$.

(c) Suppose $d>\ell+1$. Then $\t$ is injective on $\P^*$ if and
only if $\lfloor n/e\rfloor< \ell^2+\ell$.
\end{lem}

\noindent{\it Proof.} This follows from Lemmas \ref{zx1},
\ref{zx4} and \ref{zx3}.

\begin{exa}
\label{ma2} We examine the first case lying outside of the scope
of Theorem \ref{A1}, namely the case $\lfloor n/e\rfloor=d\ell+1$ and
$d\leq \ell$.

Suppose first $d<\ell$. Then $\P^*=\HP$.  The proof of Lemma
\ref{zx1} shows that $\p$ only repeats at $P=[d\ell+1,0]$ and
$Q=[0,d]$, so $V=|\P^*|-1$. The proof of Theorem \ref{A10} shows that $\phi(P)$
is the $(d+1)$th largest value of $\phi$. By Corollary \ref{icrit} all factors of (\ref{filt}) are
irreducible, except perhaps for the $(d+1)$th factor from the top,
namely $M(P)$. Since $\p$ only repeats at $P$ and $Q$, it follows
from Theorem \ref{niatro} that either $M(P)$ is irreducible, or
$M(P)=N(P)\oplus N(Q)$ with $N(P)$ and $N(Q)$ irreducible. In the
latter case $c(L)=V$ and in the former $c(L)=V-1$. Here $V=\ell d^2/2+\ell d/2+2d+1$
by Theorem \ref{A10} (else use $V=|\P^*|-1$ and Lemma \ref{longit}).
The simplest case occurs when $\ell=2$, $q=5$ and $n=6$.

Suppose next $d=\ell$. We then have
$\P^*=\HP\cup\{[0,0,1],[1,0,1]\}$. The proofs of Lemmas \ref{zx1}
and Lemma \ref{zx3} show that that $\p$ only repeats at
$P=[\ell^2+1,0,0]$ and $Q=[0,\ell,0]$, as well as
$P'=[\ell^2-\ell+1,1,0]$ and $Q'=[0,0,1]$. Thus $V=|\P^*|-2$.
The proof of Theorem \ref{A10} shows that $\phi(P')$ and $\phi(P)$ are
the $(\ell+1)$th and $(\ell+2)$th largest values of $\phi$.
The rest follows as before, except that now there are two
doubtful irreducible factors, namely $M(P')$ and $M(P)$, which are
the $(\ell+1)$th and $(\ell+2)$th factors from the top. Moreover,
$V=\ell^3/2+\ell^2/2+2\ell+2$. The simplest
example occurs when $\ell=2$, $q=3$ and $n=10$.
\end{exa}

\noindent{\bf{\Large Acknowledgments}}

I want to express my deep gratitude to Professor Dragomir Djokovic
for his invaluable contribution interpreting James' work, the
profitable time spent together discussing Gow's conjecture, and
his careful reading of the paper. I also thank Professor Rod Gow for his constant encouragement to
pursue this problem.
I am indebted to the referee for valuable suggestions. I am very thankful to Malena, Federico and Cecilia for their
continuous help and support while this research took place.

\medskip

\noindent{\bf{\Large References}}

\end{document}